\pdfoutput=1
\documentclass[11pt,final,reqno]{amsart}
\usepackage[T1]{fontenc}
\usepackage[utf8]{inputenc}
\usepackage{color}

\let\oldmarginpar\marginpar
\renewcommand\marginpar[1]{\-\oldmarginpar[\raggedleft\footnotesize #1]{\raggedright\footnotesize #1}}

\usepackage[msc-links]{amsrefs}

\newcommand{\alinea}[1]{\;#1)}

\usepackage{amsmath,amssymb}
\usepackage{mathtools}
\usepackage[all]{xy}

\theoremstyle{plain}
\newtheorem{theorem}{Theorem}[section]
\newtheorem{proposition}[theorem]{Proposition}
\newtheorem{lemma}[theorem]{Lemma}
\newtheorem{corollary}[theorem]{Corollary}
\theoremstyle{definition}
\newtheorem{definition}[theorem]{Definition}
\newtheorem{example}[theorem]{Example}
\theoremstyle{remark}
\newtheorem{remark}[theorem]{Remark}

\DeclareMathOperator{\ad}{ad}
\DeclareMathOperator{\Ad}{Ad}
\DeclareMathOperator{\AD}{\mathbf{Ad}}
\DeclareMathOperator{\id}{id}
\DeclareMathOperator{\Aut}{Aut}
\DeclareMathOperator{\tot}{tot}
\DeclareMathOperator{\pr}{pr}
\DeclareMathOperator{\ev}{ev}
\DeclareMathOperator{\Hol}{Hol}
\DeclareMathOperator{\bas}{bas}
\DeclareMathOperator{\DR}{DR}
\DeclareMathOperator{\sing}{sing}
\DeclareMathOperator{\deux}{\mathbf{2Gpd}}
\DeclareMathOperator{\boldGamma}{\boldsymbol{\Gamma}}
\DeclareMathOperator{\boldDelta}{\boldsymbol{\Delta}}
\DeclareMathOperator{\boldA}{\boldsymbol{A}}
\DeclareMathOperator{\boldB}{\boldsymbol{B}}
\DeclareMathOperator{\boldC}{\boldsymbol{C}}
\DeclareMathOperator{\boldD}{\boldsymbol{D}}
\DeclareMathOperator{\boldE}{\boldsymbol{E}}
\DeclareMathOperator{\boldF}{\boldsymbol{F}}
\DeclareMathOperator{\boldZ}{\boldsymbol{Z}}

\newcommand{\RR}{\mathbb{R}}

\newcommand{\ZZ}{\mathbb{Z}}

\newcommand{\rond}{\circ}
\newcommand{\gendex}[2]{\left\{ #1 \right\}_{#2}}
\newcommand{\genrel}[2]{\left\{ #1 | #2 \right\}}
\newcommand{\XX}{\mathfrak{X}}
\newcommand{\OO}{\Omega}

\newcommand{\mfg}{\mathfrak{g}}
\newcommand{\mfzg}{Z(\mathfrak{g})}
\newcommand{\inv}{^{-1}}
\newcommand{\toto}{\rightrightarrows}

\newcommand{\xto}[1]{\xrightarrow{#1}}
\newcommand{\xfrom}[1]{\xleftarrow{#1}}
\newcommand{\into}{\hookrightarrow}

\newcommand{\vm}{\star}
\newcommand{\hm}{*}
\newcommand{\ops}{{*}}
\newcommand{\morita}{\mathcal{M}}
\newcommand{\gm}{\Gamma}
\newcommand{\tgm}{\tilde{\Gamma}}
\newcommand{\Deltatilde}{\tilde{\Delta}}
\newcommand{\gammatilde}{\tilde{\gamma}}
\newcommand{\deltatilde}{\tilde{\delta}}
\newcommand{\dA}{d^{\boldA}}
\newcommand{\dB}{d^{\boldB}}
\newcommand{\dC}{d^{\boldC}}
\newcommand{\dD}{d^{\boldD}}
\newcommand{\partialB}{\partial^{\boldB}}
\newcommand{\partialD}{\partial^{\boldD}}
\newcommand{\lion}{p_{13}}
\newcommand{\ddR}{d_{\DR}}
\newcommand{\com}{{\scriptscriptstyle\bullet}}
\newcommand{\upcom}{^{\scriptscriptstyle\bullet}}
\newcommand{\lcom}{_{\scriptscriptstyle\bullet}}
\newcommand{\Duniv}[1]{\boldsymbol{CC}_{#1}}
\newcommand{\gpbd}{{\mathfrak{B}}}
\newcommand{\Crossed}[2]{[#1 \to #2]}
\newcommand{\DDcent}[1]{\boldsymbol{DD}_{(#1)}}
\newcommand{\ZDR}[1]{Z_{\DR}^{#1}}
\newcommand{\BDR}[1]{B_{\DR}^{#1}}
\newcommand{\GM}[5]{ \boldsymbol{#1} \underset{\sim}{\xfrom{#2}} \boldsymbol{#3} \xto{#4} \boldsymbol{#5} } 
\newcommand{\LGM}[2]{ \boldsymbol{#1} \xfrom{\sim} \boldsymbol{E}_1 \to \dots \xfrom{\sim} \boldsymbol{E}_n \to \boldsymbol{#2} } 
\newcommand{\gmto}[1]{\stackrel{#1}\rightsquigarrow} 
\begin{document}
%\pubyr{}
%\volno{}
%\iss{}

\title[Principal 2-group bundles]{$G$-gerbes, principal 2-group bundles\\ and characteristic classes}
%\date{31st century}
\author{Gr\'egory Ginot}
\address{UPMC -- Sorbonne Universit\'es -- Paris 6, 
Institut de Math\'ematiques de Jussieu--Paris Rive Gauche} 
\email{gregory.ginot@imj-prg.fr} 
\author{Mathieu Sti\'enon}
\address{Penn State University, Department of Mathematics}
\email{stienon@psu.edu}
\thanks{Research supported by the European Union through the FP6 Marie Curie R.T.N. ENIGMA
(Contract number MRTN-CT-2004-5652).} 

\begin{abstract}
Let $G$ be a Lie group and $G\to \Aut(G)$ be the canonical group homomorphism induced by the adjoint action of a group on itself.
We give an explicit description of a 1-1 correspondence between Morita equivalence classes of, on the one hand, principal 2-group $[G\to\Aut(G)]$-bundles over Lie groupoids  and, on the other hand, $G$-extensions of Lie groupoids (i.e.\ between principal $[G\to\Aut(G)]$-bundles over differentiable stacks and $G$-gerbes over differentiable
stacks). This approach also allows us to identify $G$-bound gerbes and $[Z(G)\to 1]$-group bundles over differentiable stacks, where $Z(G)$ is the center of $G$.  We also introduce universal characteristic classes for 2-group bundles. 
For groupoid central $G$-extensions, we introduce Dixmier--Douady classes that can be computed from connection-type data generalizing the ones for bundle gerbes. We prove that these classes coincide with universal characteristic classes. As a corollary, we obtain further that Dixmier--Douady classes are integral.
\end{abstract}

\maketitle
%\tableofcontents

\section{Introduction}

This paper is devoted to the relation
between groupoid $G$-extensions and \emph{principal Lie 2-group bundles} and to their characteristic classes.  

A \emph{Lie 2-group} is a Lie groupoid $\Gamma_2\toto\Gamma_1$, 
whose spaces of objects $\Gamma_1$ and of morphisms $\Gamma_2$ 
are Lie groups and all of whose structure maps are group morphisms. One shall note that in this paper we are interested in \emph{strict} Lie 2-groups only, though we believe all our results can be extended to weak ones as well.
A \emph{crossed module} $(G\xto{\rho}H)$ is a Lie group morphism $G\xto{\rho}H$ together with an action of $H$ on $G$ satisfying suitable compatibility conditions. It is standard that Lie 2-groups are in bijection with crossed modules \cites{Mackenzie, BaezLauda}. In this paper, $[G\xto{\rho} H]$ denotes the 2-group corresponding to the crossed module $(G\xto{\rho} H)$. 

Lie 2-groups arise naturally in mathematical physics. For instance, in higher gauge theory \cites{BaSc,string2group}, Lie 2-group bundles provide a well suited framework for describing the parallel transport of strings \cites{MackayPicken,AsCaJu, BaSc}. 
Several recent works have approached the concept of bundles with a ``structure Lie 2-group'' over a manifold from various perspectives \cites{Bar,BaSc,AsCaJu,SaStSc,Wockel}.
Here we take an alternative point of view and give a definition of principal Lie 2-group bundles of a global nature (i.e.\ not resorting to a description explicitly involving local charts and cocycles) and which allows for the base space to be a Lie groupoid. In other words, we consider 2-group principal bundles over differentiable stacks \cite{BehrendXu}. Our approach immediately leads to a natural construction of ``universal characteristic classes'' for principal 2-group bundles. 

Let us start with Lie (1-)groups. A principal $G$-bundle $P$ over a manifold $M$  canonically determines a homotopy class of maps from $M$ to the classifying space $BG$ of the group $G$. In fact, the set of isomorphism classes of $G$-principal bundles over $M$ is in bijection with the set of homotopy classes of maps $M\xto{f} BG$ \cites{ChernSun,Steenrod,Stasheff}. Pulling back the generators of $H^*(BG)$ (the universal classes) through $f$, one obtains characteristic classes of the principal bundle $P$ over $M$. These characteristic classes coincide with those obtained from a connection by applying the Chern--Weil construction \cites{MilnorStasheff,Dupont}.

There is an analogue but much less known, differential geometric rather than purely topological, point of view: a principal $G$-bundle over a manifold $M$ can be thought of as a ``generalized morphism'' (in the sense of Hilsum \& Skandalis \cite{Skandalis}) from the manifold $M$ to the Lie group $G$ both considered as 1-groupoids. To see this, recall that a principal $G$-bundle can be defined as a collection of transition functions $g_{ij}:U_{ij}\to G$ on the double intersections $U_{ij}$ of some open covering $\gendex{U_i}{i\in I}$ of $M$, satisfying the cocycle condition $g_{ij}g_{jk}=g_{ik}$. These transition functions constitute a morphism of groupoids from the \v{C}ech groupoid $\coprod U_{ij}\toto\coprod U_i$ associated to the open covering $\{U_i\}_{i\in I}$ to the Lie group $G\toto\ops$. Hence we have a diagram 
\[ (M\toto M)\xfrom{\sim}(\coprod U_{ij}\toto\coprod U_i)\to(G\toto\ops) \]
in the category of Lie groupoids and their morphisms whose leftward arrow is a Morita equivalence, in other words a generalized morphism from the manifold $M$ to the Lie group $G$. 

This second point of view, or more precisely its generalization to the 2-groupoid context, constitutes the foundation on which our approach is built. 
The generalization of the concept of ``generalized morphism'' to 2-groupoids
is straightforward: a generalized morphism of Lie 2-groupoids $\boldGamma \gmto{} \boldDelta$ is a diagram $\boldGamma\xleftarrow[\sim]{\phi}\boldE\xto{f}\boldDelta$ in the category $\deux$ of Lie 2-groupoids and their morphisms, where $\phi$ is a Morita equivalence (a ``smooth'' equivalence of 2-groupoids). 
It is sometimes useful to  think of two Morita equivalent Lie
2-groupoids as  two different choices of an atlas (or open cover) on
the same geometric object (which is a differentiable
2-stack~\cites{Br5, BehrendXu}).  
We define a principal $[G\xto{\rho} H]$-bundle over a Lie groupoid $\gm$ 
to be a generalized morphism from $\gm$ to $[G\xto{\rho} H]$ (up to equivalence). See Section~\ref{S:2bundle}. 

The concept of (geometric) nerve of Lie groupoids extends to the 2-categorical context as a functor from the category of Lie 2-groupoids to the category of simplicial manifolds~\cite{Street}. By convention, the cohomology of a 2-groupoid is the cohomology of its nerve, which can be computed via a double complex (for instance, see~\cite{GinotXu}). Crucially, Morita equivalences induce isomorphisms in cohomology. Therefore, any generalized morphism of 2-groupoids $\gm \gmto{F} \Crossed{G}{H}$ defining a principal $[G\to H]$-bundle $\gpbd$ over the groupoid $\gm$ yields a pullback homorphism $F^*:H\upcom(\Crossed{G}{H})\to H\upcom(\gm)$ in cohomology, which is called the \emph{cohomology characteristic map} (characteristic map for short).
The cohomology classes in $H\upcom(\Crossed{G}{H})$~\cite{GinotXu} should be viewed as universal characteristic classes and their images by $F^*$ as the characteristic classes of $\gpbd$. 
  
Lie 2-group principal bundles are closely related to  non-abelian gerbes.
Geometrically, non-abelian $G$-gerbes over differentiable stacks  can be considered as
groupoid $G$-extensions modulo Morita equivalence \cite{NADG}.
By  a groupoid $G$-extension, we mean a short exact sequence of groupoids
$1\to M\times G\xto{i}\tgm\xto{\phi}\Gamma \to 1 $,  where $M\times G$ is a bundle of groups.

One of our main results is an equivalence between (strict) $G$-gerbes over a differentiable stack and principal bundles over the 2-group $[G\to\Aut(G)]$. More precisely, we establish an explicit 1-1 correspondence between groupoid $G$-extensions up to Morita equivalence (i.e.\ strict $G$-gerbes over differentiable stacks) and principal $[G\to\Aut(G)]$-bundles over Lie groupoids modulo Morita  equivalence (i.e.\ $[G\to\Aut(G)]$-principal bundles over differentiable stacks). This is Theorem~\ref{BreenGeometric}. Note that a restricted version of this correspondence is highlighted in \cite{Friedrich}*{Theorem~4}.

It is known that Giraud's second non abelian cohomology group $H^2(\XX,G)$ classifies the $G$-gerbes over a differentiable stack $\XX$~\cite{Giraud} while Dedecker's $H^1(\XX,[G\to\Aut(G)])$ classifies the principal $[G\to\Aut(G)]$-bundles \cites{Dedecker1, Dedecker2}. 
In \cites{Breen, Br4}, Breen showed that these two cohomology groups are isomorphic. In some sense, our theorem above can be considered as an explicit geometric proof of Breen's theorem in the smooth context.
Indeed, one of the main motivations behind the present paper is the relation between $G$-extensions and 2-group principal bundles. We believe our result throws a bridge between the groupoid extension approach to the differential geometry of $G$-gerbes developed in \cite{NADG} and the one based on higher gauge theory due to Baez \& Schreiber \cite{BaSc}. This will be investigated somewhere else.

An important class of $G$-extensions is formed by the so called \emph{central $G$-extensions} \cite{NADG}, those for which the structure 2-group $\Crossed{G}{\Aut(G)}$ reduces to the 2-group $\Crossed{Z(G)}{1}$ (where $Z(G)$ stands for the center of $G$). They correspond to $G$-gerbes with trivial band or $G$-bound gerbes \cite{NADG}. Each such extension determines a principal $\Crossed{Z(G)}{1}$-bundle over the base groupoid $\gm$. 
In \cite{BehrendXu}, Behrend \& Xu gave a natural construction associating a class in $H^3(\gm)$ to a central $S^1$-extension of a Lie groupoid $\Gamma$. When the base Lie groupoid is Morita equivalent to a smooth manifold (viewed as a trivial 2-groupoid), a central $S^1$-extension is what has been studied by Murray and Hitchin under the name bundle gerbe \cites{Murray, Hitchin}. The Behrend--Xu class of a bundle gerbe coincides with its \emph{Dixmier--Douady class}, which can be described by the $3$-curvature. 
In the present paper, we extend the construction of Behrend \& Xu and define a Dixmier--Douady class $\DDcent{\alpha}\in H^3(\gm)\otimes\mfzg$ for any central $G$-extension, where $G$ is connected with a reductive Lie algebra. 
Since a central $G$-extension induces a $\Crossed{Z(G)}{1}$-principal bundle over  $\gm$, there is also a charateristic map $H^3(\Crossed{Z(G)}{1})\to H^3(\gm)$.
Dualizing, one obtains a class $\Duniv{\phi}\in H^3(\gm)\otimes\mfzg$. 
We prove that the Dixmier--Douady class $\DDcent{\alpha}$ coincides with the characteristic class $\Duniv{\phi}$. In a certain sense, this is the \emph{gerbe analogue of the Chern--Weil isomorphism} for principal bundles \cites{Dupont,MilnorStasheff}. 

The paper is organized as follows. Section~\ref{S:GMgpbd} is devoted to generalized morphisms of Lie 2-groupoids and to 2-group bundles and recalls some  standard material on  Lie 2-groupoids. The main feature of Section~\ref{S:main1} is Theorem~\ref{BreenGeometric} on the equivalence of groupoids $G$-extensions and principal $[G\xto{\Ad} \Aut(G)]$-bundles. In Section~\ref{S:CClasses} we define the characteristic map/classes of principal Lie 2-group bundles, we present the construction of the Dixmier--Douady classes of groupoid central $G$-extensions and we prove that the Dixmier--Douady class of a central $G$-extension coincides with the universal characteristic class of the induced $\Crossed{Z(G)}{1}$-bundle --- see  Theorem~\ref{T:CC=DD}. Since the universal characteristic map can be defined in cohomology with integer coefficients, we obtain that the Dixmier--Douady class of a central $G$-extension  is integral  as a corollary of our study. This applies, in particular, in the classical case of a bundle gerbe.  

Note that, when $G$ is discrete, the relation between groupoid $G$-extensions and 2-group principal bundles was also independently studied by Hae\-fli\-ger~\cite{Hae}. 

Some of the results of the present paper are related to results announced by Baez \& Stevenson~\cite{cafe}. Recently, Sati, Stasheff \& Schreiber have studied characteristic classes for 2-group bundles by the mean of $L_\infty$-algebras~\cite{SaStSc}. It would be very interesting to relate their construction to ours using integration of $L_\infty$-algebras as in~\cites{Getzler, Henriques}. 
 
\section{Generalized morphisms and principal Lie 2-group bundles}
\label{S:GMgpbd}

\subsection{Lie 2-groupoids, Crossed modules and Morita morphisms}
This section is concerned with Lie 2-groupoids and Morita equivalences. The material is rather standard. For instance, see~\cites{Mackenziebook, Moerdijkbook} for the general theory of Lie groupoids and~\cites{BaezLauda, Zhu} for Lie 2-groupoids. We only deal with the case of strict 2-groupoids and strict 2-groups.

\begin{definition}
A \emph{Lie 2-groupoid} is a double Lie groupoid 
\begin{equation}\label{double} \vcenter{\xymatrix{ \Gamma_2 \ar@<2pt>[r]^s\ar@<-2pt>[r]_t \ar@<2pt>[d]^l\ar@<-2pt>[d]_u & \Gamma_0 \ar@<2pt>[d]^\id\ar@<-2pt>[d]_\id \\ \Gamma_1 \ar@<2pt>[r]^s\ar@<-2pt>[r]_t & \Gamma_0 }} \end{equation} 
in the sense of~\cite{Mackenzie}, 
where the right column $\xymatrix{\Gamma_0 \ar@<2pt>[r]^\id\ar@<-2pt>[r]_\id & \Gamma_0}$ denotes the trivial groupoid associated to the smooth manifold $\Gamma_0$. It makes sense to use the symbols $s$ and $t$ to denote the source and target maps of the groupoid $\Gamma_2\toto \Gamma_0$ since $s\rond l=s\rond u$ and $t\rond l=t\rond u$. 
\end{definition}

\begin{remark} A Lie 2-groupoid is thus a small 2-category in which all arrows are invertible, the sets of objects, 1-arrows and 2-arrows are smooth manifolds, all structure maps are smooth and the sources and targets are surjective submersions.
\end{remark}

In the sequel, the 2-groupoid~\eqref{double} will be denoted 
$\xymatrix{ \Gamma_2 \ar@<2pt>[r]^l\ar@<-2pt>[r]_u & 
\Gamma_1 \ar@<2pt>[r]^s\ar@<-2pt>[r]_t & \Gamma_0 }$ 
or just $\boldsymbol{\Gamma}$.
The so called vertical (resp.\ horizontal) multiplication in the groupoid $\xymatrix{ \Gamma_2 \ar@<2pt>[r]^l\ar@<-2pt>[r]_u & 
\Gamma_1 }$ (resp.\ $\xymatrix{ \Gamma_2 \ar@<2pt>[r]^s\ar@<-2pt>[r]_t & 
\Gamma_0 }$) will be denoted by $\vm$ (resp.\ $\hm$)

Clearly, a Lie groupoid $\xymatrix{ \Gamma_1 \ar@<2pt>[r]^s\ar@<-2pt>[r]_t & 
\Gamma_0 }$ can be seen as a Lie 2-groupoid 
\[ \xymatrix{ \Gamma_1 \ar@<2pt>[r]^\id\ar@<-2pt>[r]_\id & 
\Gamma_1 \ar@<2pt>[r]^s\ar@<-2pt>[r]_t & \Gamma_0 } .\]

A Lie 2-groupoid where $\Gamma_0$ is the one-point space $*$ is known as a \emph{Lie 2-group}. 

There is a well-known equivalence between Lie 2-groupoids and crossed modules of groupoids~\cite{Mackenzie}. 
\begin{definition}
A \emph{crossed module} of groupoids $(X\xto{\rho}\Gamma)$ 
is a morphism of groupoids
\[ \xymatrix{X_1\ar[r]^{\rho} \ar@<0.5ex>[d]^p \ar@<-0.5ex>[d]_p & \gm_1 \ar@<0.5ex>[d]^s
\ar@<-0.5ex>[d]_t \\ X_0 \ar[r]_{\id} & \gm_0 } ,\] 
from a family of groups $X_1\xto{p}X_0$ to a groupoid 
$\Gamma_1\toto\Gamma_0$ sharing the same unit space 
$X_0=\Gamma_0$, together with a right action by automorphisms $(\gamma,x)\mapsto x^{\gamma}$ of $\Gamma_1\toto\Gamma_0$ 
on $X_1\to X_0$ satisfying:
\begin{align}
\rho(x^{\gamma})&=\gamma^{-1}\rho(x)\gamma 
&& \forall (x,\gamma)\in X_1\times_{\Gamma_0}\Gamma_1, \label{cm1} \\
x^{\rho(y)}&=y^{-1}xy && \forall (x,y)\in X_1\times_{\Gamma_0} X_1 \label{cm2}. 
\end{align}
Note that the equalities~\eqref{cm1} and~\eqref{cm2} make sense because $X_1$ is a family of groups. 
\end{definition}

\begin{example} Given any Lie group $G$, we obtain a crossed module by setting $X_1=G$, $\gm_1=\Aut(G)$, $\gm_0=\ops$ and $\rho(g)=\AD_g$ (the conjugation by $g$). 
\end{example}

\begin{example} A Lie groupoid $\gm_1\toto \gm_0$ induces a crossed module in the following way. Let $S_\gm =\{ x\in \gm_1 
| s(x)=t(x)\}$ be the set of closed loops in $\gm_1$. Then $S_\gm$ is a family of groups over $\gm_0$ and $\gm_1$ acts by conjugation on $S_\gm$. 
Therefore, we obtain a crossed module 
\[ \xymatrix{S_\gm\ar[r]^{i} \ar@<0.5ex>[d] \ar@<-0.5ex>[d] & \gm_1 \ar@<0.5ex>[d]
\ar@<-0.5ex>[d] \\ \gm_0 \ar[r]_{\id} & \gm_0 } \] 
where $i$ is the inclusion map.
\end{example}

A 2-groupoid $\xymatrix{ \Gamma_2 \ar@<2pt>[r]^l\ar@<-2pt>[r]_u & 
\Gamma_1 \ar@<2pt>[r]^s\ar@<-2pt>[r]_t & \Gamma_0}$
determines a crossed module of group\-oids $(G\xto{\rho}H)$ as follows. Here the groupoid $H$ is $\Gamma_1\toto\Gamma_0$, 
\[ G_1=\genrel{g\in\Gamma_2}{l(g)\in\Gamma_0\subset\Gamma_1},\] $\rho$ is the restriction of $u$ to $G_1$ and the action of $H_1=\gm_1$ on $G_1\subset\Gamma_2$ is by conjugation. More precisely, if $1_h$ is the unit over an object $h$ in the groupoid 
$\xymatrix{ \Gamma_2 \ar@<2pt>[r]^l\ar@<-2pt>[r]_u & \Gamma_1 }$, then   $g^h=1_{h\inv}\hm g\hm 1_h$. Conversely, given a crossed module of groupoids $(X\xto{\rho}\Gamma)$, one gets a Lie 2-groupoid 
$\xymatrix{ X_1\ltimes\Gamma_1 \ar@<2pt>[r]^l\ar@<-2pt>[r]_u & 
\Gamma_1 \ar@<2pt>[r]^s\ar@<-2pt>[r]_t & \Gamma_0 }$, where 
$X_1\ltimes\Gamma_1\toto\Gamma_1$ is the transformation groupoid and $X_1\ltimes\Gamma_1\toto\Gamma_0$ is the semi-direct product of groupoids. More precisely, for all $x,x'\in X_1$ and $\gamma,\gamma'\in\Gamma_1$, the structures maps are defined by 
\begin{align*} 
& l(x,\gamma)=\gamma, 
&& (x',\gamma')\hm(x,\gamma)=(x'x^{{\gamma'}\inv},\gamma'\gamma), \\ 
& u(x,\gamma)=\rho(x)\gamma, 
&& (x',\rho(x)\gamma)\vm(x,\gamma)=(x'x,\gamma) 
.\end{align*}

In the sequel, we will denote the Lie 2-groupoid associated to the crossed module $(G\xto{\rho} H)$ by $[G\xto{\rho} H]$.

\begin{example} The crossed module of groups $\big(G\xto{\AD}\Aut(G)\big)$ yields 
the 2-group $\xymatrix{ G\ltimes\Aut(G) \ar@<2pt>[r]^l\ar@<-2pt>[r]_u & \Aut(G) 
\ar@<2pt>[r]\ar@<-2pt>[r] & \ops }$ with structure maps
\begin{gather*}
l(g,\varphi)=\varphi \qquad u(g,\phi)=\AD_g\rond\varphi \\ 
(g_1,\AD_{g_2}\rond\varphi_2)\vm(g_2,\varphi_2)=(g_1g_2,\varphi_2) \\ 
(g_1,\varphi_1)\hm(g_2,\varphi_2)=(g_1\varphi_1(g_2),\varphi_1\rond\varphi_2)\end{gather*}
\end{example}

\smallskip

A (strict) morphism $\boldGamma\xto{\phi}\boldDelta$ of Lie 2-groupoids is a triple $(\phi_0,\phi_1,\phi_2)$ of smooth maps $\phi_i:\gm_i \to \Delta_i$ ($i=0,1,2$) commuting with all structure maps. Morphisms of crossed modules are defined similarly.

\smallskip

Let $\boldDelta$ be a Lie 2-groupoid. Given a surjective submersion  $f:M\to \Delta_0$,  we can form the \emph{pullback} Lie 2-groupoid \[ \boldDelta[M]: \,\xymatrix{ \Delta_2[M] \ar@<2pt>[r]^l\ar@<-2pt>[r]_u & 
\Delta_1[M] \ar@<2pt>[r]^s\ar@<-2pt>[r]_t &  M} ,\] where 
\[ \Delta_i[M]=\{(m,\gamma,n)\in M\times \Delta_i\times M \text{ s.t. } s(\gamma)=f(m),\; t(\gamma)=f(n) \} ,\] 
for $i\in\{1,2\}$. The maps $s,t$ are the projections on the first and last factor respectively. The maps $u,l$,
the horizontal and vertical multiplications are induced by the ones on $\boldDelta$ as follows:
\begin{align*} 
& u(m,\gamma,n)=(m,u(\gamma),n), 
&& (m,\gamma,n)\hm(n,\gamma',p)=(m,\gamma\hm\gamma',p), \\
& l(m,\gamma,n)= (m,l(\gamma),n), 
&& (m,\gamma,n)\vm(m,\gamma',n)=(m,\gamma\vm\gamma',n) 
.\end{align*}
There is a natural map of groupoids $\boldDelta[M]\to\boldDelta$ defined by $m\mapsto f(m)$ and $(m,\gamma,n)\mapsto \gamma$.

\smallskip

Pullback of 2-groupoids yield a convenient definition of  a higher analogue for Lie 2-groupoids of the notion of Morita morphism or weak equivalence of Lie 1-groupoids. These maps are higher analogues of the notion of a cover and thus are called hypercovers to agree with the terminology of~\cites{Zhu, Wolfson}.

\begin{definition}\label{D:hypercover}
A morphism of Lie 2-groupoids $ \boldGamma \xto{\phi} \boldDelta$
is a \emph{hypercover} if $\phi$ is the composition of two morphisms 
\[ \xymatrix{ 
\Gamma_2 \ar@<2pt>[d]\ar@<-2pt>[d] \ar[r] 
& \Delta_2[\Gamma_0] \ar@<2pt>[d]\ar@<-2pt>[d] \ar[r] 
& \Delta_2 \ar@<2pt>[d]\ar@<-2pt>[d] \\ 
\Gamma_1 \ar@<2pt>[d]\ar@<-2pt>[d] \ar[r] 
& \Delta_1[\Gamma_0] \ar@<2pt>[d]\ar@<-2pt>[d] \ar[r] 
& \Delta_1 \ar@<2pt>[d]\ar@<-2pt>[d] \\ 
\Gamma_0 \ar[r]^\id & \Gamma_0 \ar[r]^{\phi_0} & \Delta_0 
} \]
such that $\Gamma_0\to \Delta_0$ and 
$\Gamma_1\to\Delta_1[\Gamma_0]$ are surjective submersions and 
\[ \xymatrix{ 
\Gamma_2 \ar@<2pt>[d]\ar@<-2pt>[d] \ar[r] 
& \Delta_2[\Gamma_0] \ar@<2pt>[d]\ar@<-2pt>[d] \\ 
\Gamma_1 \ar[r] & \Delta_1[\Gamma_0] 
} \]
is a Morita morphism\footnote{Here we follow the terminology of~\cites{BehrendXu,NADG}. Moerdijk calls such maps weak equivalences~\cite{Moerdijkbook}.}of 1-groupoids.
\end{definition}

\begin{definition}
The (weakest) equivalence relation generated by the hypercovers is called \emph{Morita equivalence}.\footnote{By analogy with~\cite{BehrendXu}, we will sometimes refer to a hypercover as a Morita morphism since it induces a Morita equivalence.} 
More precisely, two Lie 2-groupoids $\boldGamma$ and $\boldDelta$ are Morita equivalent if there exists 
a finite collection $\boldE_0,\boldE_1,\dots,\boldE_n$ of Lie 2-groupoids with $\boldE_0=\boldGamma$ and $\boldE_n=\boldDelta$, and, for each $i\in\{1,\dots,n\}$, either a hypercover $\boldE_{i-1}\xto{\sim}\boldE_i$ or a hypercover $\boldE_i\xto{\sim}\boldE_{i-1}$. 
\end{definition}
In fact, by Lemma~\ref{L:MoritaClass}\alinea{\ref{seattle}}, 
one has the following well-known lemma. 

\begin{lemma}\label{L:Morita}
If $\boldGamma$ and $\boldDelta$ are Morita equivalent, there exits a chain of hypercovers $\boldGamma\xfrom{\sim}\boldE\xto{\sim}\boldDelta$ of length 2 in between $\boldGamma$ and $\boldDelta$.
\end{lemma}

\begin{remark}\label{rmk:Morita} From the categorical point of view,  
a hypercover $\phi:\boldGamma \xto{\sim} \boldDelta$ is in particular a 2-equivalence of  2-categories  preserving the smooth structures. 

We expect that the notion of Morita equivalence introduced here will help shed light on the integration problem for Courant algebroids~\cite{LWX} and more specifically on the relation between the different proposed approaches~\cites{MR2970717, MR2861783,Sheng}. It is expected that the objects integrating Courant algebroids 
are symplectic 2-stacks~\cite{arXiv:1310.6587}.
\end{remark}
 
\begin{remark}\label{rmk:stack}
Similar to~\cite{BehrendXu}, one can define differentiable 2-stacks. Two Lie 2-groupoids define the same differentiable 2-stack if, and only if, they are Morita equivalent. In fact a Lie 2-groupoid can be thought of as a choice of a differentiable atlas on a differentiable 2-stack.
\end{remark}

\subsection{Generalized morphisms of Lie 2-groupoids}\label{S:GM}
Generalized morph\-isms of Lie 2-groupoids are a straightforward generalization of generalized morphisms of Lie (1-)groupoids~\cites{Skandalis, Moerdijkbook}. They also have been considered in~\cite{Zhu}. Let $\deux$ denote the category of Lie 2-groupoids and morphisms of Lie 2-groupoids. 
\begin{definition}
A \emph{generalized morphism} $F$ is a zigzag
\[ \LGM{\Gamma}{\Delta} ,\] where all leftward arrows are hypercovers. 
We use a squig arrow $F:\boldGamma\gmto{}\boldDelta$ to denote a generalized morphism. 
The composition  of two generalized morphisms is defined by the concatenation of two zigzags. 
\end{definition}
In fact we are interested in equivalence classes of generalized morphisms:

In the sequel, we will consider two morphisms of 2-groupoids $f:\boldGamma\to \boldDelta$ and $g:\boldGamma\to \boldDelta$ to be \emph{equivalent} if there exists two smooth applications $\varphi:\gm_0\to \Delta_1$ and $\psi: \gm_1\to \Delta_2$ such that, for any $x\in\gm_2$ and any pair of composable arrows $i,j\in \gm_1$, the following relations are satisfied: 
\begin{gather} 
\label{eq:nt1}
\big( g_2(x) \hm 1_{\varphi(s(x))} \big) \vm \psi(l(x)) 
= \psi(u(x)) \vm \big( 1_{\varphi(t(x))} \hm f_2(x) \big), \\
\label{eq:nt2} 
\psi(j\hm i) = \big( 1_{g_1(j)} \hm \psi(i) \big) \vm \big( \psi(j) \hm 1_{f_1(i)} \big).
\end{gather}
In other words, $f$ and $g$ are ``conjugate'' by a (invertible) map $\psi$ compatible with the horizontal multiplication.

It is easy to check that the conditions~\eqref{eq:nt1} and~\eqref{eq:nt2} are equivalent to the data of a natural 2-transformation from $f$ to $g$~\cites{Benabou, Leinster}. Recall that a natural 2-transformation is given by the following data: an arrow $\varphi(m)\in\Delta_1$ for each object $m\in\Gamma_0$, and a 2-arrow $\psi(\gamma)\in\Delta_2$ for each arrow $\gamma\in\Gamma_1$ as in the diagram 
\[ \xymatrix{ f(s(j)) \ar[r]^{\varphi(s(j))} \ar[d]_{f(j)}& g(s(j)) \ar[d]^{g(j)}\\
f(t(j)) \ar[r]_{\varphi(t(j))} \ar@2{->}[ur]|{\psi(j)} & g(t(j)) } \]
and satisfying obvious compatibility conditions with respect to the compositions of arrows and 2-arrows. 

We now introduce the notion of equivalence of generalized morphisms;  it is the natural equivalence relation on generalized morphism extending the equivalence of groupoids morphisms. 
\begin{definition}
\emph{Equivalence of generalized morphisms} is the weakest equivalence relation satisfying the following three properties: 
\begin{enumerate}
\item If there exists a natural transformation between a pair $f,g$ of homomorphisms of 2-groupoids, $f$ and $g$ are equivalent as generalized morphisms. 
\item If $\boldGamma\xto{\phi}\boldDelta$ is a hypercover of 2-groupoids, the generalized morphisms $\boldDelta\xfrom{\phi}\boldGamma\xto{\phi}\boldDelta$ and $\boldGamma\xto{\phi}\boldDelta\xfrom{\phi}\boldGamma$ are equivalent to $\boldDelta\xto{\id}\boldDelta$ and   
$\boldGamma\xto{\id}\boldGamma$, respectively. 
\item Pre- and post-composition with a third generalized morphism preserves the equivalence.
\end{enumerate}
\end{definition}
Generalized morphisms can be seen as the $1$-morphisms in a bicategory of fractions $\deux[\morita^{-1}]$, where we have \lq\lq{}formally inverted\rq\rq{} the collection $\morita$ of hypercovers, and equivalence between generalized morphisms as being (invertible) 2-morphisms.  We refer to~\cite{Schapira}*{Chapter~7} for details on localization of categories with respect to a multiplicative system (see Lemma~\ref{L:MoritaClass} below) and to~\cite{PronkThesis}*{Chapter 3} and~\cite{Pronk96}*{Section 2} for details on the construction of the two cells of the associated bicategory of fractions. 
In particular, the $2$-morphisms in the bicategory of fractions $\deux[\morita^{-1}]$ are represented by diagrams 
\[ \xymatrix@M=6pt@R=10pt{ && {\bf F} \ar@/_2pc/@2{->}[dd]_{\alpha} \ar[dll]_{\phi_1}  \ar[drr]^{f_1} \ar@/^2pc/@2{->}[dd]^{\beta}&&    \\
\boldGamma&& {\bf E} \ar[d]^{\varepsilon_2} \ar[u]^{\varepsilon_1}  & &\boldDelta  \\ 
&& { \bf G} \ar[ull]^{\phi_2} \ar[urr]_{f_2} & &
} \] 
in which $\phi_i \circ \varepsilon_i$ are Morita morphisms (i.e.\ hypercovers) and $\alpha$, $\beta$ are 2-transformations as above. See~\cite{PronkThesis}*{Section 3.2.3} or~\cite{Pronk96}*{Section 2}.  

However, we will mosty only need the 1-category underlying $\deux[\morita^{-1}]$. 
%as constructed in Remark~\ref{R:MoritaCategory}.
\begin{example}\label{ex:GM} Let
$F_1:\boldsymbol{\gm} \underset{\sim}{\stackrel{\phi_1}{\longleftarrow}} \boldsymbol{E_1} \xto{f_1} \boldsymbol{\Delta}$ and $F_2:\boldsymbol{\gm} \underset{\sim}{\stackrel{\phi_2}{\longleftarrow}} \boldsymbol{E_2} \xto{f_2} \boldsymbol{\Delta}$ be 
two generalized morphisms. Suppose that there exists a 
morphism ${\bf E_1}\xto{\varepsilon}{\bf E_2}$ such that 
the diagram 
\[ \xymatrix@M=4pt@R=8pt{ 
& {\bf E_1} \ar[dl]_{\phi_1} \ar[dd]^{\varepsilon} \ar[dr]^{f_1} & \\ 
\boldGamma& & \boldDelta \\ 
& { \bf E_2} \ar[ul]^{\phi_2} \ar[ur]_{f_2} & 
} \] 
commutes up to 2-transformations (in particular, $\phi_2\circ \varepsilon$ is an hypercover). Then $F_1$ and $F_2$ are equivalent generalized morphisms.
\end{example}

\begin{example}\label{ex:Giso}
By its very definition, a Morita equivalence of groupoids ${\bf \gm} \xfrom{\sim} {\bf E}_1 \xto{\sim} \dots \xfrom{\sim} {\bf E}_n \xto{\sim} {\bf \Delta}$ defines two generalized morphisms $F:\boldGamma\gmto{} \boldDelta$ and $G:\boldDelta\gmto{} \boldGamma$. The compositions $F\circ G$ and $G\circ F$ are both equivalent to the identity. 
\end{example}

\begin{remark}\label{R:MoritaCategory}
As previously mentionned,  generalized morphisms are obtained by formally inverting the Morita morphisms (i.e. hypercovers). 
In fact, the following Lemma can be checked. 
\begin{lemma}\label{L:MoritaClass}
The collection $\morita$ of all hypercovers of Lie 2-groupoids is a left multiplicative 
system~\citelist{\cite{Schapira}*{Definition~7.1.5} \cite{Weibel}*{Definition~10.3.4}} in $\deux$. 
Indeed, the following properties hold:  
\begin{enumerate}
\item\label{boston} $(\boldGamma\xto{\id}\boldGamma)\in\morita$, $\forall \boldGamma\in\deux$;
\item\label{toronto} $\morita$ is closed under composition; 
\item\label{montreal} given $\boldGamma\xto{f}\boldDelta\xleftarrow[\sim]{\phi}\boldE$ in $\deux$ with $\phi\in\morita$, there exists 
$\boldGamma\xleftarrow[\sim]{\psi}\boldZ\xto{g}\boldE$ in $\deux$ with $\psi\in\morita$ such that 
\[ \xymatrix@M=4pt@R=6pt{ & \boldZ \ar[dl]_\psi^\sim \ar[dr]^g & \\ 
\boldGamma \ar[dr]_f & & \boldE \ar[dl]^\phi_\sim \\ 
& \boldDelta & } \] 
commutes; 
\item\label{seattle} given 
$\xymatrix{ \boldGamma \ar[r]^\phi_\sim & \boldDelta \ar@<2pt>[r]^f\ar@<-2pt>[r]_g & \boldE }$ 
in $\deux$ with $\phi\in\morita$, 
$f\rond\phi=g\rond\phi$ implies $f=g$.
\end{enumerate}
\end{lemma}
\begin{proof}
Properties~\alinea{\ref{toronto}} and~\alinea{\ref{montreal}} follow from~\cite{Wolfson}*{Theorem 2.12} or~\cite{Zhu}*{Section~2}. Property~\alinea{\ref{seattle}} follows from the fact that the map $\phi_0: \Gamma_0 \to \Delta_0$ is a surjective submersion and that $\phi$ is an equivalence of $2$-categories.
\end{proof}

Since $\morita$ is a left multiplicative system in the category $\deux$, we can consider the localization $\deux_\morita$ of $\deux$ with respect to $\morita$~\citelist{\cite{Schapira}*{Chapter~7} \cite{Weibel}*{Section~10.3}}. 
This new category $\deux_\morita$ has the same objects as $\deux$ but its arrows are equivalence classes 
of generalized morphisms. 
An isomorphism in $\deux_\morita$ corresponds to (the equivalence class of) a Morita equivalence in $\deux$. 
\end{remark}
In particular, the category $\deux_\morita$ is the $1$-category obtained from the bicategory $\deux[\morita^{-1}]$ (by identifying all 1-morphisms that are connected by 2-morphisms). 

Lemma~\ref{L:MoritaClass}\alinea{\ref{montreal}} implies that any generalized morphism can be represented by a chain of length 2:
\begin{lemma}\label{L:2chainGM}
Any generalized morphism between two Lie 2-groupoids $\boldGamma$ and $\boldDelta$ is equivalent to a diagram 
\[ \GM{\Gamma}{\phi}{E}{f}{\Delta} \] in the category $\deux$ in which $\phi$ is a hypercover 
(i.e.\ $\phi\in\morita$).
\end{lemma}

\begin{remark}
There is a bijection between maps of (representable) differentiable 2-stacks and equivalence classes of generalized morphisms of Lie 2-groupoids up to Morita equivalences.
\end{remark}

\subsection{Lie 2-group bundles}\label{S:2bundle}

In this section, we give a definition of Lie 2-group bundles of a global nature and formulated in terms of generalized morphisms of Lie 2-groupoids.
\begin{definition}
A \emph{principal (Lie 2-group) $[G\to H]$-bundle over a Lie groupoid $\Gamma_1\toto\Gamma_0$} is a generalized morphism 
$\gpbd$ from $\Gamma_1\toto\Gamma_0$ (seen as a Lie 2-groupoid) to the Lie 2-group $[G\to H]$ associated to the crossed module $(G\to H)$.
\end{definition}

In particular, a principal $[G\to\Aut(G)]$-bundle over a groupoid $\Gamma_1\toto\Gamma_0$ is a generalized morphism from $\Gamma_1\toto\Gamma_0$ (seen as a Lie 2-groupoid) to the Lie 2-group $[G\to\Aut(G)]$. 

Two principal $[G\to H]$-bundles $\gpbd$ and $\gpbd'$ over the Lie groupoid $\Gamma_1\toto\Gamma_0$ are said to be \emph{isomorphic} if, and only if, these two generalized morphisms are equivalent.

A $[G\to H]$-bundle over a manifold $M$ is a (2-group) $[G\to H]$-bundle over the Lie groupoid $M\toto M$.

\smallskip

Let $\gpbd$ be a  $[G\to H]$-bundle over a Lie groupoid $\gm_1\toto \gm_0$. If $\gm_1'\toto \gm_0'$ and $[G'\to H']$ are Morita equivalent to $\gm_1\toto \gm_0$ and $[G\to H]$ respectively, then the composition
\[ 
\big(\gm_1'\toto \gm_0'\big) \leftrightsquigarrow  \big(\gm_1\toto \gm_0\big) \, \gmto{\gpbd}\,[G\to H] \leftrightsquigarrow [G'\to H']
.\] 
defines a principal $[G'\to H']$-bundle over $\gm_1'\toto \gm_0'$ denoted $\gpbd$ by abuse of notation. 
Here the left and right squig arrows are the Morita equivalences seen as invertible generalized morphisms as in Example~\ref{ex:Giso}. 
\begin{definition}\label{D:2groupbundle}
A principal (Lie 2-group) $[G\to H]$-bundle $\gpbd$ over a Lie groupoid $\Gamma_1\toto\Gamma_0$ and a principal (Lie 2-group) $[G'\to H']$-bundle $\gpbd'$ over a Lie groupoid $\Gamma'_1\toto\Gamma'_0$ are said to be \emph{Morita equivalent} if, and only if,  
$\Gamma'_1\toto\Gamma'_0$ is Morita equivalent to $\Gamma_1\toto\Gamma_0$, $[G'\to H']$ is Morita equivalent to $[G\to H]$, 
and $\gpbd$ (viewed as a generalized morphism $\big(\gm_1'\toto \gm_0'\big)\gmto{} [G'\to H']$) and $\gpbd'$ are equivalent generalized morphisms. 
\end{definition} 

\begin{remark}
When the groupoid is just a manifold, our definition is equivalent to the usual definition of Lie 2-group bundles in~\cites{BaSc,Bar,SaStSc,Wockel} as suggested by Examples~\ref{E:principal} and~\ref{ex:MGAUTG} below.
\\
Furthermore our notion of isomorphism of principal 2-group bundles over a fixed manifold agrees with the one of~\cite{Wockel} and is in fact a particular case of the one described in~\cite{Wolfson}, see Remark~\ref{R:Wolfson}.
\end{remark}

\begin{example}\label{E:principal}
Let $P\xto{\pi} M$ be a principal $H$-bundle. Then the diagram \[ \xymatrix{ M \ar@<2pt>[d]\ar@<-2pt>[d] 
& P\times_M  P\ar[l] \ar@<2pt>[d]\ar@<-2pt>[d] \ar[r]
& H \ar@<2pt>[d]\ar@<-2pt>[d] \\ 
M \ar@<2pt>[d]\ar@<-2pt>[d] 
& P\times_M P\ar[l]_{\phi} \ar@<2pt>[d]^s\ar@<-2pt>[d]_t \ar[r]^f 
& H \ar@<2pt>[d]\ar@<-2pt>[d] \\ 
M & P \ar[l]_\pi \ar[r] & \ops 
} \] where $s(x,y)=x$, $t(x,y)=y$, $\pi(x)=\phi(x,y)=\pi(y)$ and $x\cdot f(x,y)=y$, 
defines a generalized morphism from the manifold $M$ to the $2$-group $\Crossed{1}{H}$. Hence, it is a $2$-group bundle over $M$. 
Note that a principal $H$-bundle $P$ over $M$ is Morita equivalent (as a 2-group bundle) to a principal $H'$-bundle $P'$ over $M'$ if, and only if, $H$ and $M$ are isomorphic to $H'$ and $M'$ respectively and $P$ and $P'$ are isomorphic principal bundles.
\end{example}

\begin{example}\label{ex:MGAUTG}
Let $M$ be a smooth manifold and $G$ be a (non-abelian) Lie group. 
A non abelian 2-cocycle~\cites{Giraud,Moerdijk,Dedecker1,Dedecker2} on $M$ with values in $G$ relative to an open covering $\gendex{U_i}{i\in I}$ of $M$ is a collection of smooth maps 
\[ \lambda_{ij}: U_{ij}\to\Aut(G) \qquad \text{and} \qquad 
g_{ijk}:U_{ijk}\to G \] 
satisfying the following relations: 
\begin{gather*}
\lambda_{ij}\rond\lambda_{jk}=\AD_{g_{ijk}}\rond\lambda_{ik} \\ 
g_{ijl}g_{jkl}=g_{ikl}\lambda_{kl}\inv(g_{ijk}) 
.\end{gather*}
Such a non-abelian 2-cocycle defines a $[G\to\Aut(G)]$-bundle over the manifold $M$; for it can be seen as the generalized morphism 
\[ \xymatrix{ 
M \ar@<2pt>[d]\ar@<-2pt>[d] 
& \coprod_{i,j}U_{ij}\times G\times G \ar[l] \ar@<2pt>[d]^l\ar@<-2pt>[d]_u \ar[r]^f 
& G\ltimes\Aut(G) \ar@<2pt>[d]\ar@<-2pt>[d] \\ 
M \ar@<2pt>[d]\ar@<-2pt>[d]  
& \coprod_{i,j}U_{ij}\times G \ar[l]_{\phi} \ar@<2pt>[d]\ar@<-2pt>[d] \ar[r] 
& \Aut(G) \ar@<2pt>[d]\ar@<-2pt>[d] \\ 
M & \coprod_{i}U_{i} \ar[l] \ar[r] & \ops 
} \] 
between the manifold $M$ and the 2-group $[G\to\Aut(G)]$.
Here 
\begin{align*}
& l(x_{ij},g_1,g_2)=(x_{ij},g_1) && \phi(x_{ij},g)=x \\
& u(x_{ij},g_1,g_2)=(x_{ij},g_2) && f(x_{ij},g_1,g_2)=\big(g_2 g_1\inv,\AD_{g_1}\rond\lambda_{ij}(x)\big)
\end{align*}
where $x_{ij}$ denotes a point $x\in M$ seen as a point of the open subset $U_{ij}=U_i\cap U_j$, $x_{i}$ the point $x\in M$ seen as a point of the open subset $U_i$, and $g,g_1,g_2$ arbitrary elements of $G$. 
The horizontal and vertical multiplication are given by
\begin{align*} &(x_{ij},g_1) \hm (x_{jk},g_2) = \big(x_{ik}, g_{ijk}\lambda_{jk}^{-1}(g_1)g_2\big),\\
&(x_{ij},g_1,g_2)\vm (x_{ij},g_2,g_3) = (x_{ij},g_1,g_3). \end{align*}
\end{example}

\begin{example}
Let $\gendex{U_i}{i\in I}$ be an open covering of a smooth manifold $M$. A family of smooth maps $g_{ijk}:U_{ijk}\to S^1$ defines a Lie groupoid structure on $\coprod_{i,j}U_{ij}\times S^1\toto\coprod_{i}U_{i}$ 
with multiplication 
\[ (x_{ij},e^{i\varphi})\cdot(x_{jk},e^{i\psi})=
(x_{ik},g_{ijk}e^{i(\varphi+\psi)}) \]  
if, and only if, $g_{ijk}$ is a \v{C}ech 2-cocycle.
In that case, we get the generalized morphism of 2-groupoids 
\[ \xymatrix{ 
M \ar@<2pt>[d]\ar@<-2pt>[d] 
& \coprod_{i,j}U_{ij}\times S^1\times S^1 \ar[l] \ar@<2pt>[d]^l\ar@<-2pt>[d]_u \ar[r]^{\qquad \quad f} 
& S^1 \ar@<2pt>[d]\ar@<-2pt>[d] \\ 
M \ar@<2pt>[d]\ar@<-2pt>[d]  
& \coprod_{i,j}U_{ij}\times S^1 \ar[l] \ar@<2pt>[d]\ar@<-2pt>[d] \ar[r] 
& \ops \ar@<2pt>[d]\ar@<-2pt>[d] \\ 
M & \coprod_{i}U_{i} \ar[l] \ar[r] & \ops 
} \] 
with $f(x_{ij},e^{i\varphi},e^{i\psi})=e^{i(\psi-\varphi)}$.
It defines an $[S^1\to\ops]$-bundle over $M$.
\end{example}
 
\begin{remark}\label{R:Wolfson}
Our definition of principal 2-groups bundles over a stack also agrees with one which was  introduced more recently by Wolfson in~\cite{Wolfson}; more precisely our definition agrees with a special case of principal bundle over a simplicial Lie group 
in \emph{loc.~cit.}\ when the simplicial Lie group is a strict Lie 2-group (viewed as a special kind of simplicial Lie groups as in~\cite{Wolfson}*{\S 6}). Indeed, given a generalized morphism $ \big(\gm_1\toto \gm_0\big) \, \gmto{\gpbd}\,[G\to H]$, we get a strict Lie 2-groupoid map $\boldsymbol{U} \to [G\to H]$, for some hypercover 
$\boldsymbol{U}{\stackrel{\phi}\longrightarrow}\boldsymbol{\Gamma}$, 
which in turns yields a map between the associated simplicial manifolds. Pulling back  the universal bundle $W([G\to H])$ 
of~\cite{Wolfson}*{Definition 6.2} along this map yields a twisted Cartesian product $\boldsymbol{U}\times_{[G\to H]} W([G\to H])\to\boldsymbol{U}$ over (the simplicial manifold associated to) $\boldsymbol{U}$ (which is a stack by construction). 
Since $\boldsymbol{U}{\stackrel{\phi}\longrightarrow}\boldsymbol{\Gamma} $ is a hypercover, this is precisely the data of a  local $2$-bundle which is principal with respect to (the simplicial Lie group associated to) $[G\to H]$ in the sense of~\cite{Wolfson}*{\S 5}. 
 
Furthermore, by choosing a common refinement of two hypercovers, an isomorphism of principal 2-group bundles in the sense of Definition~\ref{D:2groupbundle} yields an equivalence of local 2-bundles (as twisted cartesian product and as stacks), 
since in both cases it boils down to being a collection of  equivalences of local objects of the form 
$\boldsymbol{E}\cong\boldsymbol{U}\times [G\to H]$.
\end{remark}

\begin{example}\label{ex:inertia}
Principal $2$-group bundles arise whenever one studies group actions on stacks, which, in general are only \emph{weak} actions. 
For instance, there is a canonical  (but subtle) weak action of $S^1$ on the inertia stack of any differentiable stack giving rise 
to a canonical principal bundle, which shall be detailed elsewhere. 
For the moment we just explain briefly how it can be defined in terms of Lie $2$-group(oid)s.
Let $\boldsymbol{\gm}:\gm_1\toto\gm_0$ be a Lie groupoid. 
Its inertia groupoid\footnote{The inertia groupoid of a groupoid $\Gamma_1\toto\Gamma_0$ is a groupoid representing the 
inertia stack of the quotient stack $[\Gamma_0/\Gamma_1]$.  
See~\cite{BGNX} for details on inertia groupoids and stacks.} 
$\Lambda\boldsymbol{\gm}:S_\gm\times_{\gm_0}\gm_1\toto S_{\gm}$ 
(where $S_{\gm}=\{\gamma \in \gm_1 \text{ s.t.\ } d(\gamma)=s(\gamma)\}$ is the space of loops) 
has a \emph{canonical} action of the group stack associated to the 2-group $[\mathbb{Z}\to 1]$. This action is given, 
for $(\gamma,g) \in S_{\gm}\times \Gamma_1$ and $n\in \mathbb{Z}$ by $(\gamma, g)\cdot n:= (\gamma, \gamma^n\cdot g)$. 
Hence it also inherits an action of the group stack $S^1 \cong [\mathbb{Z}\to\mathbb{R}]$ (induced by the canonical map $\mathbb{R}\to 0$). One shall note that this action is almost never represented by a strict action of the group $S^1$ 
on the inertia groupoid but really by an action of the $2$-group $[\mathbb{Z}\to\mathbb{R}]$.
Assume the inertia groupoid is a Lie groupoid --- which is true if $\boldsymbol{\gm}$ is \'etale and proper. 
It follows from~\cite{GinotNoohi}*{Theorem 0.2} that the quotient of the stack represented by $\Lambda\boldsymbol{\gm}$ 
by the (group stack represented by) $S^1\cong [\mathbb{Z}\to \mathbb{R}]$ is a differentiable stack and further 
a $[\mathbb{Z} \to \mathbb{R}]$-principal bundle.
Indeed this quotient stack can be presented by a Lie groupoid $\boldsymbol{\tilde{\gm}}$  which is Morita equivalent 
to the Lie 2-groupoid $\boldsymbol{\tilde{\tilde{\gm}}}: S_\gm \times_{\gm_0} \gm_1\times \mathbb{Z}\times \mathbb{R} 
\toto S_\gm \times_{\gm_0} \gm_1\times \mathbb{R} \toto S_\gm $. 
The horizontal multiplications are given by the product of the Lie (2-)groupoids structures of $\boldsymbol{\gm}$ and 
$[\mathbb{Z}\to \mathbb{R}]$ while vertical mutiplication is induced by the action of $\mathbb{Z}$ on the inertia groupoid described above. : 
The canonical projection  $\boldsymbol{\tilde{\tilde{\gm}}} \to [\mathbb{Z} \to \mathbb{R}]= \mathbb{Z}\times \mathbb{R} \toto\mathbb{R} \toto 1$ gives right to the generalized morphism 
$\mathcal{B}:\boldsymbol{\tilde{\gm}}\xleftarrow{\sim}\boldsymbol{\tilde{\tilde{\gm}}} \to [\mathbb{Z} \to \mathbb{R}]$ 
 defining the bundle structure. 
\end{example}

\begin{remark}\label{topologicalside}
There is a nerve functor from Lie 2-groupoids to simplicial spaces generalizing the nerve for Lie 1-groupoids. 
For instance, see~\cites{Street, Segal, BuFaBl} and Section~\ref{S:cohomology} below. 
Composing it with the (fat) realization functor, we obtain the classifying space functor $\boldGamma\mapsto B\gm$ 
from Lie 2-groupoids to topological spaces. Since the realization of a Morita morphism (i.e.\ hypercover) 
is a homotopy equivalence, a generalized morphism $\boldGamma \gmto{F} \boldDelta$ induces a map 
$B\boldGamma \xto{BF} B\boldDelta$ in the homotopy category of topological spaces.
In particular, a $\Crossed{G}{H}$-group bundle over a manifold induces a  map $M \to B\Crossed{G}{H}$ in the homotopy category. 
This is the topological side of generalized morphisms and 2-group bundles. In fact, using standard arguments on homotopy for manifolds, it should be possible to prove that $\Crossed{G}{H}$-group bundles over $\boldGamma$ (up to Morita equivalences) are in bijection with homotopy classes of maps $B\boldGamma \to B\Crossed{G}{H}$. 
\end{remark}

\section{Groupoid $G$-extensions}
\label{S:main1}
We fix a Lie group $G$. We recall the following definition (see~\cite{NADG})
\begin{definition}
A \emph{Lie groupoid $G$-extension} is a short exact sequence 
of Lie groupoids over the identity map on the unit space $M$ 
\begin{equation}\label{extension}
1\to M\times G\xto{i}\tgm\xto{\phi}\Gamma \to 1 
\end{equation}
Here both $\Gamma$ and $\tgm$ are Lie groupoids over $M$ and $M\times G\toto M$ is a (trivial) bundle of groups.
\end{definition}
The map $\phi$ being a map over the identity map on the unit space $M$ means that both $\tgm$ and $\Gamma$ 
have $M$ for unit space and that the restriction of $\phi$ to $M$ is the identity on $M$: 
\[ \xymatrix{ \tgm \ar@<2pt>[d]\ar@<-2pt>[d] \ar[r]^{\phi} 
& \Gamma \ar@<2pt>[d]\ar@<-2pt>[d] \\ M \ar[r]_\id & M .} \]

In the sequel, an extension like~\eqref{extension} will be denoted $\tgm\xto{\phi}\Gamma\toto M$ and we will write $g_m$ instead of $i(m,g)$.

Lie groupoid $G$-extensions can be interpreted in terms of crossed modules as follows.
\begin{proposition}\label{vancouver} The morphism of groupoids $\tgm\xto{\phi}\Gamma\toto M$ is a groupoid $G$-extension if, and only if, $(M\times G\xhookrightarrow{i}\tgm)$ is a crossed module of groupoids with quotient groupoid $\tgm/{i(M\times G)}$ isomorphic to $\Gamma$.
\end{proposition}
Proposition~\ref{vancouver} follows easily from Remark~\ref{R:GtoGM} and Lemma~\ref{2to1} below.

\begin{definition}[\cite{NADG}]
A \emph{Morita morphism between Lie groupoid $G$-ex\-ten\-sions} is a homomorphism of Lie groupoid $G$-extensions 
\[ \xymatrix{ 
\tgm \ar[d]_f \ar[r] & \Gamma \ar@<2pt>[r]\ar@<-2pt>[r] \ar[d]_{f} & M \ar[d]_{f} \\ \Deltatilde \ar[r] & \Delta \ar@<2pt>[r]\ar@<-2pt>[r] & N 
} \]
such that $M\xto{f}N$ is a surjective submersion and 
\[ \vcenter{\xymatrix{ 
\Gamma \ar@<2pt>[d]\ar@<-2pt>[d] \ar[r]^{f} & \Delta \ar@<2pt>[d]\ar@<-2pt>[d] \\ 
M \ar[r]^{f} & N 
}} \qquad \text{and} \qquad 
\vcenter{\xymatrix{ 
\tgm \ar@<2pt>[d]\ar@<-2pt>[d] \ar[r]^{f} & \Deltatilde \ar@<2pt>[d]\ar@<-2pt>[d] \\ 
M \ar[r]^{f} & N 
}} \]
are Morita morphisms\footnote{Weak equivalences in~\cite{Moerdijkbook} and hypercovers in~\cites{Zhu,Wolfson}.} of 1-groupoids. 
\end{definition}

As in the Lie 2-groupoid case, the Morita morphisms of Lie groupoid extensions form a left multiplicative system in the category of Lie groupoid extensions and homomorphisms of Lie groupoid extensions. 
Hence, one can localize this category by its Morita morphisms. Two Lie groupoid extensions are \emph{Morita equivalent} if they are isomorphic in the localized category. 
(As in the Lie 2-groupoids case, there is a notion of \emph{generalized morphisms} for Lie groupoid extensions. In that language, a Morita equivalence is an invertible generalized morphism.) 

Here is our first main theorem.

\begin{theorem}\label{BreenGeometric}
There exists a bijection between the Morita equivalence class\-es of Lie groupoid $G$-extensions
and the Morita equivalence classes of $[G\to\Aut(G)]$-bundles over Lie groupoids. 
\end{theorem}

\begin{remark}
The above theorem can be regarded as a geometric version of a theorem of Breen~\cite{Breen}, 
which states that $H^2(\XX,G)$ is isomorphic to $H^1(\XX,(G\to\Aut(G)))$. 
\end{remark}

The proof of Theorem~\ref{BreenGeometric} is the object of the next two sections.

\subsection{From groupoid $G$-extensions to $[G\to\Aut(G)]$-bundles} 
\label{S:Gextensiontogpbd}
Given a Lie groupoid G-extension $\xymatrix{ \tgm \ar[r]^{\phi} & 
\Gamma \ar@<2pt>[r]^a\ar@<-2pt>[r]_b & M }$, one can define a Lie 2-groupoid 
$\xymatrix{ \tgm\times_\Gamma\tgm \ar@<2pt>[r]^l\ar@<-2pt>[r]_u & 
\tgm \ar@<2pt>[r]^s\ar@<-2pt>[r]_t & M }$, where 
\begin{gather*}
\tgm\times_\Gamma\tgm = \genrel{(\gammatilde_1,\gammatilde_2)\in
\tgm\times\tgm}{\phi(\gammatilde_1)=\phi(\gammatilde_2)} \\ 
l(\gammatilde_1,\gammatilde_2)=\gammatilde_2 \qquad u(\gammatilde_1,\gammatilde_2)=\gammatilde_1 \\ 
s(\gammatilde)=a\big(\phi(\gammatilde)\big) \qquad 
t(\gammatilde)=b\big(\phi(\gammatilde)\big) \\ 
(\gammatilde_1,\gammatilde_2)\vm(\gammatilde_2,\gammatilde_3) 
=(\gammatilde_1,\gammatilde_3) \\ 
(\gammatilde_1,\gammatilde_2)\hm(\deltatilde_1,\deltatilde_2) 
=(\gammatilde_1\cdot\deltatilde_1,\gammatilde_2\cdot\deltatilde_2) 
.\end{gather*}
Here $\cdot$ stands for the multiplication in $\tgm\toto M$.

The groupoid homomorphism $\phi$ naturally induces a Morita morphism (i.e.\ hypercover) of 2-groupoids: 
\begin{equation}\label{eq:MoritaExt} 
\vcenter{\xymatrix{ \tgm\times_\Gamma\tgm \ar@<2pt>[d]^l\ar@<-2pt>[d]_u \ar[r]
& \Gamma \ar@<2pt>[d]^{\id}\ar@<-2pt>[d]_{\id} \\ 
\tgm \ar@<2pt>[d]^s\ar@<-2pt>[d]_t \ar[r]^{\phi} & \Gamma \ar@<2pt>[d]^a\ar@<-2pt>[d]_b \\ 
M \ar[r]^{\id} & M }} \end{equation}
where the 2-groupoid 
\[ \xymatrix{ \Gamma \ar@<2pt>[r]^{\id}\ar@<-2pt>[r]_{\id} & 
\Gamma \ar@<2pt>[r]^a\ar@<-2pt>[r]_b & M } \] is simply the 1-groupoid 
$\xymatrix{ \Gamma \ar@<2pt>[r]^a\ar@<-2pt>[r]_b & M }$ seen as a 2-groupoid in the trivial way.

Consider the map $\tgm\to\Aut(G):\gammatilde\mapsto\AD_{\gammatilde}$ 
defined by $\big(\AD_{\gammatilde}g\big)_{t(\gammatilde)}= 
\gammatilde\cdot g_{s(\gammatilde)}\cdot\gammatilde\inv$.
It gives a morphism of Lie groupoids 
\begin{equation}\label{star} \vcenter{\xymatrix{ \tgm \ar[r]^{\AD} \ar@<2pt>[d]\ar@<-2pt>[d] & 
\Aut(G) \ar@<2pt>[d]\ar@<-2pt>[d] \\ M \ar[r] & \ops }} \end{equation}
which, together with the map 
\[ \tgm\times_\Gamma\tgm\to G\ltimes\Aut(G) : 
(\gammatilde_1,\gammatilde_2)\mapsto(g,\AD_{\gammatilde_2}), \] 
where $\gammatilde_1\gammatilde_2\inv= g_{t(\gammatilde_1)}$, defines 
a homomorphism of Lie 2-groupoids 
\begin{equation}\label{eq:TildeG} \vcenter{\xymatrix{ 
\tgm\times_\Gamma\tgm \ar@<2pt>[d]\ar@<-2pt>[d] \ar[r] & 
G\ltimes\Aut(G) \ar@<2pt>[d]\ar@<-2pt>[d] \\ 
\tgm \ar@<2pt>[d]\ar@<-2pt>[d] \ar[r]^{\AD} & \Aut(G) \ar@<2pt>[d]\ar@<-2pt>[d] \\ 
M \ar[r] & \ops 
}} \end{equation}

\begin{remark} Note that the induced map 
\[ \xymatrix{ \tgm\times_\Gamma\tgm \ar@<2pt>[d]\ar@<-2pt>[d] \ar[r] & G\ltimes\Aut(G) \ar@<2pt>[d]\ar@<-2pt>[d] \\ \tgm \ar[r] & \Aut(G) } \] 
is a fully faithful functor. 
\end{remark}

\begin{remark}\label{R:GtoGM} In terms of crossed modules, the above discussion goes as follows.
The extension $1\to M\times G\xto{i}\tgm\xto{\phi}\Gamma \to 1 $ leads to an action of $\tgm$ on the groupoid $M\times G\toto M$ by conjugation, i.e.\ via the map $\gammatilde\mapsto\AD_{\gammatilde}$. Then  $\xymatrix{ \tgm\times_\Gamma\tgm \ar@<2pt>[r]^l\ar@<-2pt>[r]_u & 
\tgm \ar@<2pt>[r]^s\ar@<-2pt>[r]_t & M }$ is the Lie 2-groupoid corresponding to the crossed module $(M\times G \xto{i} \tgm)$. The 
projection onto the first factor $M\times G\to M$ and the morphism $\phi:\tgm\to\Gamma$ induce the Morita equivalence of crossed modules $(M\times G\to \tgm) \to (M\to \Gamma)$ corresponding to the map~\eqref{eq:MoritaExt}. Moreover, the map $\AD: \tgm\to \Aut(G)$ yields the map of crossed modules $\big(M\times G \xto{i} \tgm \big) \xto{(\pr_2,\AD)} \big(G\to \Aut(G)\big)$ corresponding to the morphism of Lie 2-groupoids~\eqref{eq:TildeG}.
\end{remark} 

\begin{proposition}\label{P:extensiontoGM} 
\begin{enumerate}
\item\label{detroit} A Lie groupoid $G$-extension 
$\tgm\to\Gamma\toto M$ 
induces a principal $[G\to\Aut(G)]$-bundle over $\Gamma\toto M$, 
which can be described explicitely by the following generalized morphism: 
\[ \xymatrix{ 
\Gamma \ar@<2pt>[d]\ar@<-2pt>[d] 
& \tgm\times_\Gamma\tgm \ar[l] \ar@<2pt>[d]\ar@<-2pt>[d] \ar[r] 
& G\ltimes\Aut(G) \ar@<2pt>[d]\ar@<-2pt>[d] \\ 
\Gamma \ar@<2pt>[d]\ar@<-2pt>[d]  
& \tgm \ar[l] \ar@<2pt>[d]\ar@<-2pt>[d] \ar[r] 
& \Aut(G) \ar@<2pt>[d]\ar@<-2pt>[d] \\ 
M & M \ar[l] \ar[r] & \ops 
} \] 
\item\label{chicago} If $\tgm\to\Gamma\toto M$ and $\Deltatilde\to\Delta\toto N$ are Morita equivalent $G$-extensions, then the corresponding $2$-group bundles are Morita equivalent.
\end{enumerate}
\end{proposition}

\begin{proof}
Claim\alinea{\ref{detroit}} follows from the above discussion. Suppose given a Morita morphism of $G$-extensions \[ \xymatrix@R=12pt{ 
\tgm \ar[d]_f \ar[r] & \Gamma \ar@<2pt>[r]\ar@<-2pt>[r] \ar[d]_{f} & M \ar[d]_{f} \\ \Deltatilde \ar[r] & \Delta \ar@<2pt>[r]\ar@<-2pt>[r] & N .
} \] Since $f$ commutes with the $\tgm$ and $\Deltatilde$-actions on $G$, there is a commutative diagram
\[ \xymatrix{ [M\to\Gamma] \ar[d]_{(f,f)}^{\sim} & [M\times G\to\tgm] 
\ar[d]_{(f\times \id,f)}^{\sim} \ar[l]_{\sim} \ar[r]^{(\pr_2,\AD)} & [G\to \Aut(G)] \ar[d]^{\id}\\
[N\to\Delta] & [N\times G\to\Deltatilde] \ar[l]_{\sim} \ar[r]_{(\pr_2,\AD)} & [G\to \Aut(G)] } ,\]
where $p_2$ denotes the canonical projection on the second component. 
Now, Claim\alinea{\ref{chicago}} follows from Example~\ref{ex:GM}.
\end{proof}

\subsection{From $[G\to\Aut(G)]$-bundles to groupoid $G$-extensions}
\label{S:gpbdtoGextension}

In this section, we show how to reverse the procedure. Starting from a $[G\to\Aut(G)]$-bundle, we recover 
a groupoid $G$-extension. 

For future reference, we state the following technical result without proof. 

\begin{lemma}\label{2to1}
Let \[ \xymatrix{ \Delta_2 \ar@<2pt>[d]\ar@<-2pt>[d] \ar[r]^{\phi_2} & \Gamma_2 \ar@<2pt>[d]\ar@<-2pt>[d] \\ 
\Delta_1 \ar@<2pt>[d]\ar@<-2pt>[d] \ar[r]^{\phi_1} & \Gamma_1 \ar@<2pt>[d]\ar@<-2pt>[d] \\ 
\Delta_0 \ar[r]^{\phi_0} & \Gamma_0 } \]
be a hypercover of Lie 2-groupoids. And let 
\[ \xymatrix{ L \ar[d]_j \ar[r]^\phi & K \ar[d]^i \\ \Delta_1 \ar[r]_{\phi_1} & \Gamma_1 } \] be the induced map of crossed modules. 
Then $\phi$ maps $j\inv(1_m)$ onto $i\inv(1_{\phi(m)})$ bijectively (for every $m\in\Delta_0$) and induces a functor 
from the groupoid $\tfrac{\Delta_1}{j(L)}$ to the groupoid $\tfrac{\Gamma_1}{i(K)}$, which is fully faithful  and surjective 
on the objects.\footnote{The crossed modules $[1\to \tfrac{\Delta_1}{j(L)}]$ and $[1\to \tfrac{\Gamma_1}{i(K)}]$ are the `cokernels' of the crossed modules $[L\xto{j} \Delta_1]$ and $[K\xto{i} \gm_1]$ respectively.}
\end{lemma}

Now, given a $[G\to\Aut(G)]$-bundle over a Lie groupoid $\Gamma\toto\Gamma_0$, 
we proceed with the construction of a  Lie groupoid $G$-extension. 
Suppose the $[G\to\Aut(G)]$-bundle is given by the generalized morphism of 2-groupoids 
\[ \xymatrix{ 
\Gamma \ar@<2pt>[d]\ar@<-2pt>[d] 
& \Delta_2 \ar[l]_{\phi_2} \ar@<2pt>[d]\ar@<-2pt>[d] \ar[r]^{f_2}
& G\ltimes\Aut(G) \ar@<2pt>[d]\ar@<-2pt>[d] \\ 
\Gamma \ar@<2pt>[d]\ar@<-2pt>[d]  
& \Delta_1 \ar[l]_{\phi_1} \ar@<2pt>[d]\ar@<-2pt>[d] \ar[r]^{f_1} 
& \Aut(G) \ar@<2pt>[d]\ar@<-2pt>[d] \\ 
\Gamma_0 & \Delta_0 \ar[l]_{\phi_0} \ar[r]^{f_0} & \ops 
} \] 
and let \[ \xymatrix{ \Gamma_0 \ar[d]^i & L \ar[l]_{\phi} \ar[d]^j \ar[r]^f & G \ar[d]^{\AD} \\ 
\Gamma & \Delta_1 \ar[l]^{\phi_1} \ar[r]_{f_1} & \Aut(G) } \]
be the induced generalized morphism of crossed modules. 
Hence \[ L=\{\alpha\in\Delta_2 | u(\alpha)=1_x \text{ for some } x\in\Delta_0\} \] 
and $j:L\to\Delta_1$ is the restriction of the structure map $l:\Delta_2\to\Delta_1$ to $L$. 
Since $\phi$ is a hypercover and $\Gamma_0\xto{i}\Gamma$ is an injection, by Lemma~\ref{2to1}, $L\xto{j}\Delta_1$ is also injective and 
\[ \xymatrix{ \tfrac{\Delta_1}{j(L)} \ar@<2pt>[d]\ar@<-2pt>[d] 
\ar[r]^{\phi} & \tfrac{\Gamma}{i(\Gamma_0)} \ar@<2pt>[d]\ar@<-2pt>[d] \\ 
\Delta_0 \ar[r]_{\phi} & \Gamma_0 } \] 
is a fully faithful functor. 
Since $\frac{\Gamma}{i(\Gamma_0)}$ is diffeomorphic to $\Gamma$, the groupoid $\frac{\Delta_1}{j(L)}\toto\Delta_0$ 
is the pullback of $\Gamma\toto\Gamma_0$ through the surjective submersion $\Delta_0\xto{\phi}\Gamma_0$. 
Therefore $\frac{\Delta_1}{j(L)}=\Delta_0\times_{\phi,\Gamma_0,s}\Gamma\times_{t,\Gamma_0,\phi}\Delta_0$ is a 
smooth manifold. 

Consider the groupoid structure on $\Delta_1\times G\toto\Delta_0$ with source $s(\delta,g)=s(\delta)$, 
target $t(\delta,g)=t(\delta)$, mutliplication 
\[ (\delta_1,g_1)\cdot(\delta_2,g_2)=\big(\delta_1\delta_2,f_1(\delta_2^{-1})[g_1]\cdot g_2\big) ,\] 
and inverse $(\delta,g)^{-1}=\big(\delta^{-1},f_1(\delta)[g^{-1}]\big)$ 
for all $\delta,\delta_1,\delta_2\in\Delta_1$ and $g,g_1,g_2\in G$. 

The map $\mathcal{H}:L\to\Delta_1\times G$ defined by $\mathcal{H}(\alpha)=\big(j(\alpha),f(\alpha^{-1})\big)$ 
is a morphism of groupoids from $L\toto\Delta_0$ to $\Delta_1\times G\toto\Delta_0$. 
One checks that 
\begin{align*} 
(\delta,g)\cdot\big(j(\alpha),f(\alpha^{-1})\big) 
&=\big(j(\delta * \alpha * \delta^{-1}),f(\delta * \alpha * \delta^{-1})^{-1}\big)\cdot(\delta,g) 
\\ &=\big(\delta\cdot j(\alpha),f(\alpha^{-1})\cdot g\big) 
,\end{align*} 
for all $\delta\in\Delta_1$, $g\in G$, and $\alpha\in L$. 
Thus the image of $L$ under $\mathcal{H}$ is a normal subgroupoid of $\Delta_1\times G\toto\Delta_0$. 
Since $j:L\to\Delta_1$ is injective, the action of $L$ on $\Delta_1\times G$ by multiplication from the right 
\[ (\delta,g)\bullet\alpha=(\delta,g)\cdot\big(j(\alpha),f(\alpha^{-1})\big)=\big(\delta\cdot j(\alpha),f(\alpha^{-1})\cdot g\big) \] 
is free and its orbit space $(\Delta_1\times G)/\mathcal{H}(L)$ is a smooth manifold. 

Now consider the groupoid $G$-extension 
$\Delta_1\times G\to\Delta_1\toto\Delta_0$. 
The morphism of groupoids $\Delta_1\times G\ni(\delta,g)\mapsto\delta\in\Delta_1$ intertwines the right action of $L$ 
on $\Delta_1\times G$ with the right action $\delta\bullet\alpha=\delta\cdot j(\alpha)$ of $L$ on $\Delta_1$, 
whose orbit space is the smooth manifold $\frac{\Delta_1}{j(L)}$.
Therefore, passing to quotients, we obtain the $G$-extension of groupoids 
\[ (\Delta_1\times G)/\mathcal{H}(L)\to \Delta_1/j(L)\toto\Delta_0 .\] 
Note that the corresponding crossed module is \[ \big(\Delta_0\times G\to (\Delta_1\times G)/\mathcal{H}(L)\big) .\]

\begin{proposition}\label{P:gpbdtoGextension}
\begin{enumerate} 
\item\label{Cleveland} Every $[G\to\Aut(G)]$-bundle 
over a Lie groupoid $\Gamma\toto\Gamma_0$ induces a Lie groupoid $G$-extension. 
\item\label{Columbus} Morita equivalent $[G\to\Aut(G)]$-bundles induce Morita equivalent extensions. 
\end{enumerate}
\end{proposition}

\begin{proof} 
\alinea{\ref{Cleveland}} As was outlined above, a bundle 
\[ [\Gamma_0\to\Gamma]\xfrom{\phi}\boldDelta\xto{f}[G\to\Aut(G)] \] 
determines a $G$-extension 
\[ (\Delta_1\times G)/\mathcal{H}(L)\to\Delta_1/j(L)\toto\Delta_0 ,\] 
where $j:L\to\Delta_1$ is the restriction of the structure map $l:\Delta_2\to\Delta_1$ to 
\[ L=\{\alpha\in\Delta_2 | u(\alpha)=1_x \text{ for some } x\in\Delta_0\} \] and 
$\mathcal{H}:L\to\Delta_1\times G$ is the morphism of groupoids from $L\toto\Delta_0$ to $\Delta_1\times G\toto\Delta_0$
defined by $\mathcal{H}(\alpha)=\big(j(\alpha),f(\alpha^{-1})\big)$.

\alinea{\ref{Columbus}} It is sufficient to check that  for any diagram
\[ \xymatrix@M=4pt@R=8pt{ 
& \boldE \ar[dl]_{\phi_1}^{\sim} \ar[dd]^{\varepsilon} \ar[dr]^{f_1} & \\ 
[\Gamma_0\to\Gamma] & & [G\xto{\AD} \Aut(G)] \\ 
& \boldF \ar[ul]^{\phi_2}_{\sim} \ar[ur]_{f_2} & 
} \] commuting up to natural 2-equivalences, the $G$-extensions corresponding to the lower and upper generalized morphisms are Morita equivalent. Since $\phi_1$, $\phi_2$ are Morita equivalences, $\varepsilon$ is also a Morita equivalence. Therefore, by Lemma~\ref{L:2chainGM}, we can assume that $\varepsilon$ is a hypercover.
Then, denoting by $(K\to E_1)$ and $(L\to F_1)$ the crossed modules corresponding to $\boldE$ and $\boldF$ respectively, the map $\varepsilon$ induces a commutative diagram
\[ \xymatrix{ (E_1\times G)/\mathcal{H}(K) \ar[d]_{(\varepsilon,\id)} \ar[r] & E_1/j(K) \ar@<2pt>[r]\ar@<-2pt>[r] \ar[d]_{\varepsilon} 
& E_0 \ar[d]_{\varepsilon} \\ (F_1\times G)/\mathcal{H}(L) \ar[r] & F_1/j(L) \ar@<2pt>[r]\ar@<-2pt>[r] & F_0 } \]
which is a Morita equivalence of extensions by Lemma~\ref{2to1}. 
\end{proof}

\subsection{Proof of Theorem~\ref{BreenGeometric}}
It remains to prove that the constructions of Section~\ref{S:Gextensiontogpbd} and Section~\ref{S:gpbdtoGextension} are inverse of each other.

Suppose that a $[G\to \Aut(G)]$-principal bundle $\gpbd$ over ${\gm}$ is given by the generalized morphism 
$\gm\xfrom{\phi}\boldDelta\xto{f}[G\to\Aut(G)]$. 
Let $\frac{\Delta_1\times G}{\mathcal{H}(L)}\to\frac{\Delta_1}{j(L)}\toto\Delta_0$ be the induced $G$-principal extension as in Proposition~\ref{P:gpbdtoGextension}\alinea{\ref{Cleveland}}. 
The corresponding crossed module is $(\Delta_0\times G\to \frac{\Delta_1\times G}{\mathcal{H}(L)})$. 
We have the following commutative diagram of crossed modules: 
\[ \xymatrix{ (\Gamma_0\to\Gamma) & (L\to\Delta_1) \ar[l]_{\sim} \ar[d] \ar[r] & (G\to\Aut(G)) 
\\ (\Delta_0\to\frac{\Delta_1}{j(L)}) \ar[u]^{\sim} & (\Delta_0\times G\to\frac{\Delta_1\times G}{\mathcal{H}(L)}) \ar[l]^{\sim} \ar[ru] & } ,\]
where $(L\to\Delta_1)$ is the crossed module corresponding to $\boldDelta$. 
It follows that the generalized morphism 
\[ [\Gamma_0\to\Gamma]\xfrom{\sim}[L\to\Delta_1]\to[G\to\Aut(G)] \] we started from is  equivalent to the generalized morphism 
\[ [\Gamma_0\to\Gamma]\xfrom{\sim}[\Delta_0\times G\to\tfrac{\Delta_1\times G}{\mathcal{H}(L)}]\to [G\to\Aut(G)] \] associated to the $G$-extension $\frac{\Delta_1\times G}{\mathcal{H}(L)}\to\frac{\Delta_1}{j(L)}\toto\Delta_0$.
Hence they represent the same $(G\to\Aut(G))$-bundle over $\Gamma\toto\Gamma_0$.

Reciprocally, if $\tgm\to \gm \toto M$ is a $G$-extension, then the associated principal $[G\to \Aut(G)]$-bundle is given by the generalized morphism \begin{equation}\label{eq:GtoGMtoG}[M\to \gm] \xfrom{\phi} [M\times G \to \tgm] \xto{(\pr_2,\AD)} [G\to \Aut(G)]\end{equation} according to Remark~\ref{R:GtoGM}. Direct inspection of the proof of Proposition~\ref{P:gpbdtoGextension} shows that the $G$-extension induced 
by the generalized morphism~\eqref{eq:GtoGMtoG} is exactly $\tgm\to \gm \toto M$. 

\section{Universal characteristic maps and Dixmier--Douady classes}
\label{S:CClasses}
\subsection{Cohomology of Lie 2-groupoids}
\label{S:cohomology}
To each Lie 2-groupoid $\Gamma_2\toto\Gamma_1\toto\Gamma_0$ is associated a simplicial manifold: 
its \emph{(geometric) nerve} $N_\com\boldGamma$.  It is the nerve of the underlying  2-category as defined by Street~\cite{Street}.  
In particular, $N_0\boldGamma=\gm_0$, $N_1\boldGamma=\gm_1$, $N_2\boldGamma$ is a submanifold of $\gm_2\times \gm_1\times \gm_1\times \gm_1$ parameterizing the 2-arrows of $\gm_2$ fitting in a commutative triangle
\begin{equation}\label{eq:N2}\vcenter{\xymatrix{ & A_1\ar[rd]^{f_0} \ar@{=>}[d]^{\alpha}  &  \\ 
A_0\ar[rr]_{f_1} \ar[ru]^{f_2} &&A_2 }}
\end{equation}
and $N_3\boldGamma$ is a submanifold of $(\gm_2)^{4}\times (\gm_1)^{6}$ parameterizing the commutative tetrahedra like 
\begin{equation}\label{eq:tetrahedron} 
\vcenter{\xymatrix@C=3pt@R=12pt@M=6pt{ &&& A_3 &&& \\ &&&& && \\ &&& A_1 \ar[uu]^{f_{02}} \ar@{=>}[ul]^{\alpha_2} \ar@{=>}[d]^{\alpha_3} \ar[drrr]^{f_{03}} & \ar@{=>}[ul]_{\alpha_0} && \\ A_0 \ar[rrrrrr]_{f_{13}} \ar[urrr]^(0.65){f_{23}} \ar[rrruuu]^{f_{12}} && \ar@{:>}[ul]_<{\alpha_{1}} &&&& A_2 \ar[llluuu]_{f_{01}} }}
\end{equation} 
with faces given by elements of $N_2{\gm}$. 
By the commutativity of the tetrahedron~\eqref{eq:tetrahedron}, we mean that $(\alpha_3\hm f_{01})\vm\alpha_1=(f_{23}\hm\alpha_0)\vm\alpha_2$. 
For $p\geq 3$, $N_p{\bf \gm}$ is the manifold of all $p$-simplices such that each subsimplex of dimension $3$ is a tetrahedron of the form~\eqref{eq:tetrahedron} above \cites{Street,Duskin}.  
The nerve of a Lie groupoid considered as a Lie 2-groupoid is isomorphic to its usual (1-)nerve \cite{Segal}. 
The nerve $N_\com$ defines a functor from the category of Lie $2$-groupoids to the category of simplicial manifolds. Taking the fat realization of the nerve defines a functor from Lie $2$-groupoids to topological spaces. 

\smallskip

The de~Rham cohomology groups of a Lie 2-groupoid $\boldGamma$ are defined to be the total cohomology groups of the bicomplex $(\Omega^\com(N_\com\boldGamma),\ddR,\partial)$, where $\ddR:\Omega^p(N_q\boldGamma)\to\Omega^{p+1}(N_q\boldGamma)$ stands for the de~Rham differential and $\partial:\Omega^p(N_q\boldGamma)\to
\Omega^p(N_{q+1}\boldGamma)$ is defined by $\partial=(-1)^p\sum_{i=0}^{q+1}(-1)^i d_i^*$, where $d_i:N_{q+1}\boldGamma\to N_{q}\boldGamma$ denotes the $i$\textsuperscript{th} face map. 
We use the shorter notation $\Omega^{\com}_{\tot}(\boldGamma)$ for the associated total complex. 
Hence $\Omega^{n}_{\tot}(\boldGamma)=\bigoplus_{p+q=n}\OO^p(N_q\boldGamma)$ with (total) differential $\ddR+\partial$. 
We denote the subspaces of cocycles and coboundaries by $\ZDR\upcom(\boldGamma)$ and $\BDR\upcom(\boldGamma)$ respectively, and the cohomology of $\boldGamma$ by $H^\com(\boldGamma)$.

The following Lemma is folklore (see~\cite{DuIs} for a more general statement with respect to hypercovers).
\begin{lemma}\label{lm:homGM}
Let $F:\boldGamma\to\boldDelta$ be a hypercover of Lie 2-groupoids. Then $F^*:H\upcom(\boldDelta)\to H\upcom(\boldGamma)$ is an isomorphism.
\end{lemma}
 
\begin{proof}
It is well-known that a natural transformation between two 2-functors $f$ and $g$ from $\boldGamma$ to $\boldDelta$ induces a simplicial homotopy between $f_*: N\lcom(\boldGamma) \to N\lcom(\boldDelta)$ and $g_*: N\lcom(\boldGamma) \to N\lcom(\boldDelta)$, for instance see~\cite{BuFaBl}*{Proposition 4}. In particular equivalent (topological) 2-categories have homotopic nerves. The result follows for a hypercover with a section. Since local sections always exist, the   general case  reduces to  a hypercover  $\boldGamma[\coprod U_i] \to \boldGamma$ induced by pullback along the canonical map $\coprod U_i \to \boldGamma$ where $(U_i)$ is a cover of $\boldGamma_0$. The result follows from  a classical Mayer-Vietoris argument as in~\cite{Behrend/cohomology}.
\end{proof}

By Lemma~\ref{lm:homGM} above, a generalized morphism $F:{\boldGamma} \underset{\sim}{\xfrom{\phi_1}} {\bf E}_1 \xto{f_1} \dots \underset{\sim}{\xfrom{\phi_n}} {\bf E}_n \xto{f_n} {\boldDelta} $ induces a pullback map in cohomology
\[ F^*: H\upcom(\boldDelta)\xto{f_n^*}H\upcom(\boldE_n) \xto{(\phi_n^*)^{-1}}\dots \xto{f_1^*} H\upcom(\boldE_1) \xto{(\phi_1^*)^{-1}} H\upcom(\boldGamma) \] 
Clearly, $(F\circ G)^* =G^* \circ F^*$ and, if $F$ is a Morita equivalence, then $F^*$ is an isomorphism.  

\begin{lemma}
If $F$ and $G$ are equivalent generalized morphisms from $\boldGamma$ to $\boldDelta$, the maps $F^*$ and $G^*$, which they induce at the cohomology level, are equal.
\end{lemma}

\begin{proof}  A natural transformation between two 2-functors $f$ and $g$ from $\boldGamma$ to $\boldDelta$ induces a simplicial homotopy between $f_*: N\lcom(\boldGamma) \to N\lcom(\boldDelta)$ and $g_*: N\lcom(\boldGamma) \to N\lcom(\boldDelta)$ (see~\cite{BuFaBl}). Therefore the lemma follows from the definition of equivalence of generalized morphisms and Lemma~\ref{lm:homGM}. 
\end{proof} 

\smallskip

\begin{remark}\label{paramN2}
Note that, for a Lie 2-groupoid $\boldGamma: \vcenter{\xymatrix{ \Gamma_2 \ar@<2pt>[r]^l\ar@<-2pt>[r]_u & \Gamma_1 \ar@<2pt>[r]^s\ar@<-2pt>[r]_t & \Gamma_0 }}$, $N_2\boldGamma$ may be identified to $\Gamma_2\times_{s,\Gamma_0,t}\Gamma_1$ 
so that the face maps take the form \begin{gather*}
d_0:N_2\boldGamma\to N_1\boldGamma:(\alpha,c)\mapsto u(\alpha) \\ 
d_1:N_2\boldGamma\to N_1\boldGamma:(\alpha,c)\mapsto l(\alpha)\cdot c \\ 
d_2:N_2\boldGamma\to N_1\boldGamma:(\alpha,c)\mapsto c 
\end{gather*}
More precisely, 

\[ \boxed{\xymatrix{ A_2 & A_1 
\ar@/_2pc/[l]_{u(\alpha)}="1"
\ar@/^2pc/[l]^{l(\alpha)}="2"
& A_0 \ar[l]_{c} 
 \ar@{}"1";"2"|(.2){\,}="7"
 \ar@{}"1";"2"|(.8){\,}="8"
 \ar@{=>}"7" ;"8"^{\alpha}
}} \in \Gamma_2\times_{s,\Gamma_0,t}\Gamma_1 \] 
is identified to 
\[ \boxed{\vcenter{\xymatrix{ & A_1\ar[ld]_{u(\alpha)} \ar@{=>}[d]^{\alpha * c}  &  \\ 
A_2 & & A_0 \ar[ll]^{l(\alpha) * c} \ar[lu]_{c} }}}
\in N_2\boldGamma .\] 
\end{remark}

\begin{remark}\label{rm:coefficient}
The singular cohomology $H^\com_{\sing}(\boldGamma,R)$ of $\boldGamma$ with coefficient in a ring $R$ is defined similarly to the de Rham cohomology. More precisely it is the cohomology of 
$(C_{\sing}^\com(N_\com\boldGamma),d_{\sing}+\partial)$ where $(C_{\sing}^\com(X,R),d_{\sing})$ denotes the singular cochain complex of a space $X$ with coefficients in $R$. As for manifolds, for $R=\mathbb{R}$, one has  a natural isomorphism $H^\com_{\sing}(\boldGamma,\mathbb{R})\cong H^\com(\boldGamma)$. 
\end{remark}

\subsection{Cohomological characteristic map for 2-group bundles}
\label{S:characteristicmap}

Fix a crossed module $(G\to H)$ and let $\gpbd$ be a principal $[G\to H]$-bundle over $\boldGamma$. In this section we construct a universal characteristic homomorphism $\Duniv{\gpbd}: H\upcom(\Crossed{G}{H}) \to H\upcom({\boldGamma})$ generalizing the usual characteristic classes of a principal bundle.

By definition, $\gpbd$ is a generalized morphism $\boldGamma \gmto{} \Crossed{G}{H}$. 
Therefore, passing to cohomology, we obtain the homomorphism
\begin{equation}\label{eq:Duniv}
\Duniv{\gpbd}: H\upcom(\Crossed{G}{H}) \xto{\gpbd^*}  H\upcom(\boldGamma)
\end{equation} 
which we call \emph{the characteristic homomorphism} of the $\Crossed{G}{H}$-bundle $\gpbd$. It depends only on the isomorphism class of the $2$-group bundle.
\begin{proposition}\label{P:Duniv} If $\gpbd$ and $\gpbd'$ are isomorphic $[G\to H]$-bundles over $\boldGamma$, then $\Duniv{\gpbd}=\Duniv{\gpbd'}:H\upcom(\Crossed{G}{H})\to H\upcom(\boldGamma)$.
\end{proposition}

\begin{proof} It is an immediate consequence of Lemma~\ref{lm:homGM} since isomorphic principal 2-group bundles 
are equivalent as generalized morphisms.
\end{proof} 

\begin{remark} By analogy with the case of principal bundles, one can think of the elements of $H\upcom(\Crossed{G}{H})$ 
as universal characteristic classes and their images in $H\upcom(\boldGamma)$ by $\Duniv{\gpbd}$ as characteristic classes 
of the $[G\to H]$-bundle over $\boldGamma$. 

For instance, it is proved~\cite{GinotXu}*{Proposition 6.3} that the characteristic classes associated to the string $2$-group associated to a compact simple Lie group coincide with the usual ones modulo the Pontryagin class.
\end{remark}

\begin{example}
Let $P\xto{\pi} M$ be a principal $H$-bundle. Then, by Example~\ref{E:principal}, 
$P$ induces a structure of $[1\to H]$-bundle over $M$.
Since $H\upcom(\Crossed{1}{H})\cong H\upcom(BH)$, the characteristic map $\Duniv{P}$ of this bundle coincides
with the classical map $H\upcom(BH)\to H\upcom(M)$ induced by the principal $H$-bundle 
 structure on $P$. In particular, for a compact Lie group $H$,  the characteristic map coincides with the Chern-Weil map 
 $S(\mathfrak{h}^*)^\mathfrak{h} \to H\upcom(M)$ induced by the choice of a connection on $P$. 
\end{example}

\begin{example}
From example~\ref{ex:inertia}, we know that the inertia groupoid of a Lie groupoid $\boldsymbol{\gm}$ gives rise to a 
principal $[\mathbb{Z}\to \mathbb{R}]$-bundle.
In that case $H^\bullet([\mathbb{Z}\to \mathbb{R}])\cong H^\bullet(BS^1)\cong 
\mathbb{R}[x]$ where $x$
is a generator of degree $2$. In particular, we get a characteristic class $\Duniv{\Gamma}(x) \in H^2([\Lambda \boldsymbol{\gm}/S^1])= 
H^2_{S^1}(\Lambda\boldsymbol{\gm})$ (see~\cite{GinotNoohi} for equivariant cohomology of stacks).
 For instance if $\boldsymbol{\gm}$ is the groupoid $G\toto 1$ with $G$ a simply connected compact Lie group, then its inertia groupoid is 
 the transformation Lie groupoid $G\times G\toto G$ with $G$ acting on itself by the adjoint action. From the Gysin sequence in
 equivariant homology of stacks (\cite{GinotNoohi}*{\S 8}) and the fact that the homology $H\upcom(G\times G \toto G)$ is 
 trivial in degrees $1$ and $2$, we see that  $\Duniv{G\times G \toto G}(x)$ is an integral\footnote{More precisely, 
it is the image of a generator of the cohomology with integer coefficients. See Remark~\ref{rm:CCcoefficient}.} generator of 
 $H^2([\Lambda\boldsymbol{\gm}/S^1])$. 
 \end{example}
 
Let $\boldGamma\gmto{F}\boldDelta$ be a generalized morphism  of Lie (1-)groupoids and let $\gpbd:\boldDelta \xleftarrow[\sim]{\phi}\boldE\xto{f}\Crossed{G}{H}$ be a $2$-group bundle with base $\boldDelta$. The pullback $F^*(\gpbd)$ of the $\Crossed{G}{H}$-bundle $\gpbd$ from $\boldDelta$ to $\boldGamma$ by $F$ is the composition $\gpbd\circ F$ of the two generalized morphisms. It is a principal $\Crossed{G}{H}$-bundle over $\boldGamma$.

The Whitney sum of two 2-group bundles is defined as follows. 
Let $\gpbd:\boldGamma\xleftarrow[\sim]{\phi}\boldE\xto{f}\Crossed{G}{H}$ and 
$\gpbd':\boldGamma\xleftarrow[\sim]{\phi'}\boldE'\xto{f'}\Crossed{G'}{H'}$ be two $2$-group bundles over the same base $\boldGamma$. 
Let $\boldF$ be the ``fiber product'' $2$-groupoid 
$E_2\times_{\gm_2}E_2' \toto E_1\times_{\gm_1}E_1'\toto E_0\times_{\gm_0} E_0'$ 
with the obvious structure maps: $s(e,e') =(s(e),s(e'))$, $(x,x')\hm (y,y')=(x\hm x', y\hm y')$, etc.\ 
The Whitney sum $\gpbd\oplus\gpbd'$ is the $[G\times G'\to H\times H']$-bundle over $\boldGamma$ 
given by the generalized morphism
$\boldGamma\underset{\sim}{\xfrom{\phi = \phi'}}\boldF\xto{f\times f'}\Crossed{G\times G'}{H\times H'}$.

By Proposition~\ref{P:Duniv}, we obtain
\begin{corollary}
\begin{enumerate}
\item $\Duniv{F^*(\gpbd)} = F^* \circ \Duniv{\gpbd}$.
\item $\Duniv{\gpbd\oplus \gpbd'}= \Delta^*\circ \big(\Duniv{\gpbd} \times \Duniv{\gpbd'} \big)$, 
where $\Delta:\boldGamma \to \boldGamma \times \boldGamma$ is the diagonal map 
and $\times$ is the cross-product $H\upcom(\Crossed{G}{H})\otimes H\upcom(\Crossed{G'}{H'})
\cong H\upcom(\Crossed{G\times G'}{H\times H'})$. 
\end{enumerate}
\end{corollary}

\begin{remark}\label{rm:CCcoefficient}
The result of this section easily extends to singular cohomology with any coefficient (see Remark~\ref{rm:coefficient}). 
In particular the characteristic map \[ \Duniv{\gpbd}: H\upcom_{\rm sing}(\Crossed{G}{H},\mathbb{Z}) \xto{\gpbd^*}  H\upcom(\boldGamma,\mathbb{Z}) \] is defined in cohomology with integer coefficient.
\end{remark}

\begin{remark}
By Proposition~\ref{P:extensiontoGM}, a Lie groupoid $G$-extension $\tgm\xto{\phi}\Gamma\toto M$ 
induces a principal $[G\to\Aut(G)]$-bundle $\gpbd_\phi$ over the groupoid $\Gamma\toto M$.
Hence we obtain the universal characteristic map $\Duniv{\gpbd_\phi}:H\upcom([G\xto{\AD}\Aut(G)])\to H\upcom(\boldGamma)$.
Unfortunately, the cohomology $H\upcom([G\xto{\AD}\Aut(G)])$ is not known when the center of $G$ is large and it is trivial
when the center of $G$ is of dimension less than three~\cite{GinotXu}. 
Therefore one cannot have much hope of getting interesting characteristic classes except for extensions 
whose structure 2-group can be reduced. Indeed, this is the object of the next section.
\end{remark}

\subsection{DD classes for groupoid central $G$-extensions}\label{S:DDcent}

Let $\tgm\xto{\phi} \gm \toto M$ be a $G$-extension of Lie groupoids.
Let $\phi'$ denote the factorization of 
the morphism $\phi$ through the projection $q:\tgm\to\tgm/Z(G)$: 
\[ \xymatrix{ \tgm \ar[r]^{\phi} \ar[d]_q & \Gamma \\ \tgm/Z(G) \ar[ru]_{\phi'} & } \]
The extension $\phi$ is said to be \emph{central} \cite{NADG}
if there exists a section 
$\sigma:\Gamma\to\tgm/Z(G)$ 
of $\phi'$ such that 
\begin{equation}\label{central} x g=g x 
\qquad \forall x\in q\inv\big(\sigma(\Gamma)\big), \;
\forall g\in G\end{equation}
In this case, the subspace $\tgm'=q\inv\big(\sigma(\Gamma)\big)$ of $\tgm$ 
is a central $Z(G)$-extension of $\Gamma\toto M$.
\begin{remark} The definition of a central $G$-gerbe here is taken from~\cite{NADG}. 
According to \emph{loc.~cit}, it agrees with the one of a $G$-gerbe with \emph{trivial band}. 
\end{remark}

Given $\gamma\in\tgm$, there exists $x\in\tgm'$ such that $\phi(x)=\phi(\gamma)$. 
Thus there exists $k\in G$ such that $\gamma=x\cdot k$. Given $\gamma$, 
both $x$ and $k$ are uniquely determined up to an element of $Z(G)$. 
Defining a homomorphism of Lie groupoids $r:\tgm\to G/Z(G)$ by the relation 
$q(\gamma)=\sigma(\phi(\gamma))r(\gamma)$, we obtain that, for any $g\in G$, 
\[ g\gamma=g x k=x g k=x k \cdot k\inv g k=\gamma g^{r(\gamma)} \] 
where $g^{r(\gamma)}$ denotes the conjugate $k\inv g k$ of $g$ 
by any element $k\in G$ such that $k Z(G)=r(\gamma)$.

\begin{proposition}\label{P:red}
Let $\tgm\xto{\phi}\Gamma\toto M$ be a $G$-extension of a Lie groupoid $\Gamma$ 
and let $\mathfrak{B}$ denote the corresponding $[G\to\Aut(G)]$-bundle. 
The extension is central if, and only if, the $[G\to\Aut(G)]$-bundle $\mathfrak{B}$ 
reduces to a principal $[Z(G)\to 1]$-bundle, i.e.\ there exists a generalized morphism 
$Z\mathfrak{B}:[M\to\Gamma]\to [Z(G)\to 1]$ such that 
\[ \xymatrix{ [M\to\Gamma] \ar[r]^{\hspace{-0.8cm}\mathfrak{B}} \ar[rd]_{Z\mathfrak{B}} & 
[G\to\Aut(G)] \\ & [Z(G)\to 1] \ar@{^{(}->}[u] } \] 
is commutative up to equivalence.
\end{proposition}
In particular, being central is invariant under Morita equivalences of Lie groupoids extension.
\begin{proof}
Let $\tgm\xto{\phi}\Gamma\toto M$ be a central $G$-extension. 
The corresponding 2-group bundle $\gpbd$ is the generalized morphism 
$[M\to \gm] \xfrom{} [M\times G\xto{i} \tgm] \xto{} [G\xto{\Ad} \Aut(G)]$, 
see Proposition~\ref{P:extensiontoGM} and Remark~\ref{R:GtoGM}. 
Let $\tau:\tgm'\to \tgm$ be the inclusion map. 
The $Z(G)$-extension defines the crossed module $[M\times Z(G)\xto{i'} \tgm' ]$ 
and we have a commutative diagram
\begin{equation}\label{eq:red}\vcenter{\xymatrix{[M\to \gm] 
& [M\times G\xto{i} \tgm]\ar[l] \ar[r] & [G\xto{\Ad} \Aut(G)] \\
&[M\times Z(G)\xto{i'} \tgm' ]\ar[lu] \ar[u]_{\tau} \ar[r]& [Z(G) \to 1] \ar@{^{(}->}[u] }} \end{equation} 
Note that the right square in~\eqref{eq:red} is commutative because the extension is central. 
Diagram~\eqref{eq:red} implies that the 2-group bundle $\gpbd$ reduces. 

Reciprocally, assume $\gpbd$ reduces. By Proposition~\ref{P:gpbdtoGextension}\alinea{\ref{Columbus}}, 
passing to a Morita equivalent groupoid, we can assume that the $G$-extension is the extension corresponding 
to the generalized morphism $[M\to\gm]\gmto{Z\gpbd} [Z(G)\to 1] \hookrightarrow [G\to \Aut(G)]$. 
If $Z\gpbd$ is the generalized morphism $[M\to \gm] \xfrom{} [M\times L \to \Delta] \xto{} [Z(G)\to 1]$, 
the associated extension is, according to Section~\ref{S:gpbdtoGextension}, 
$\tgm \to \gm \toto M$, where $\tgm=\big(\Delta \times_{L} Z(G)\big)\times_{Z(G)} G$. 
Since the composition $\Delta\to 1\to\Aut(G)$ is trivial, $\AD_{\widetilde{\gamma}}$ is trivial 
for all $\widetilde{\gamma}\in\tgm$. Therefore, the extension is central.
\end{proof}

\begin{remark}\label{reductive}
If $G$ is a  Lie group whose Lie algebra $\mfg$ is reductive, its Lie algebra  decomposes as a direct sum 
$\mfg\cong \mfzg \oplus \mathfrak{m}$ of ideals, where $\mfzg$ is the center of $\mfg$.
In the sequel, the symbol $\pr$ will denote the induced projection $\mfg\to\mfzg$, 
which is a homomorphism of Lie algebras and maps $[\mfg,\mfg]$ onto 0. 
Moreover, if $G$ is connected, this direct sum decomposition is not only $\ad_{\mfzg}$-invariant but also $\Ad_G$-invariant and, 
consequently, $\pr\rond\Ad_g=\pr$ for all $g\in G$. 
Moreover, for any $g\in G$ and any smooth path $t\mapsto f_t$ in $G$ with $f_0=1$ 
and $\tfrac{d}{dt}f_t\big|_0=\xi\in\mfg$, one has 
\begin{equation}\label{iets} \pr\Big( \tfrac{d}{dt}f_t\inv g f_t g\inv\big|_0\Big) 
= \pr(\Ad_g\xi-\xi)=\pr(\xi)-\pr(\xi)=0 \end{equation}
\end{remark}

\begin{proposition}\label{P:DDcent} 
Let $\tgm\xto{\phi}\Gamma\toto M$ be a central $G$-extension with $G$ connected and whose Lie algebra is
reductive.\footnote{Such Lie groups are called reductive, though this terminology sometimes applies only to algebraic groups. Examples of Lie groups with reductive Lie algebras include $GL_n(\RR)$ and all compact Lie groups.} Let $\alpha\in\Omega^1(\tgm;\mfg)$ be a connection 1-form for the right principal
$G$-bundle $\tgm\xto{\phi}\Gamma$. 
\begin{enumerate}
\item\label{BatonRouge} Then there exists $\Omega_\alpha\in\ZDR{3}(\Gamma\lcom;\mfzg)$ such that 
$\pr\big(d\alpha +\partial\alpha\big)=\phi^*(\Omega_\alpha)$.
\item\label{NewOrleans} Moreover, if $\alpha_1$ and $\alpha_2$ are two different connection 1-forms, then $\Omega_{\alpha_1}-\Omega_{\alpha_2}\in \BDR{3}(\Gamma\lcom;\mfzg)$. \end{enumerate}
We call $\DDcent{\alpha}:= [\Omega_{\alpha}] \in H^3(\Gamma)\otimes \mfzg$ the Dixmier--Douady class of the $G$-central extension.
\end{proposition}

\begin{proof}
\alinea{\ref{BatonRouge}} Being a connection 1-form, $\alpha\in\Omega^1(\tgm;\mfg)$ enjoys the following two properties: 
\begin{gather*} 
R_g^*\alpha=\Ad_{g\inv}\rond\alpha,\quad\forall g\in G \\ 
\alpha(\hat{\xi}_x)=\xi,\quad\forall x\in M,\forall\xi\in\mfg 
\end{gather*}
Given any $\xi\in\mfg$ and any $G$-invariant vector field $v\in\XX(\tgm)$, we get 
\begin{multline*} d\alpha(\hat{\xi},v)=\hat{\xi}\big(\alpha(v)\big)
-v\big(\alpha(\hat{\xi})\big)-\alpha\big([\hat{\xi},v]\big) 
=\hat{\xi}\big(\alpha(v)\big)-v(\xi)-\alpha(\mathcal{L}_{\hat{\xi}}v) \\ 
= \mathcal{L}_{\hat{\xi}}\big(\alpha(v)\big)=-\ad_{\xi}\big(\alpha(v)\big)
\end{multline*}
since the vector field $v$ is $G$-invariant and the function $\alpha(v)$ is $G$-equivariant. It follows that $\pr\rond d\alpha(\hat{\xi},v)=\pr [\alpha(v),\xi] =0$ since $\pr [\mfg,\mfg]=0$.
Moreover, we have 
\[ R_g^*(d\alpha)=d(R_g^*\alpha)=d(\Ad_{g\inv}\rond\alpha)=\Ad_{g\inv}\rond d\alpha \] for all $g\in G$.
Therefore, by Remark~\ref{reductive}, the 2-form $\pr\rond d\alpha\in\OO^2(\tgm,\mfzg)$ is basic; there exists $\omega\in\OO^2(\gm,\mfzg)$ such that $\pr\rond d\alpha=\phi^*\omega$.

Consider \[ \tgm_2=\tgm\times_{s,\gm,t}\tgm=\genrel{(x,y)\in\tgm\times\tgm}{s(x)=t(y)} ,\] the three face maps 
\[ p_1(x,y)=x \qquad m(x,y)=x\cdot y \qquad p_2(x,y)=y \] 
from $\tgm_2$ to $\tgm$ and the action of $G\times G$ on $\tgm_2$ given by 
\[ (x,y)^{(g,h)}=(x g_{s(x)},y h_{s(y)}) .\]
Then we have 
\[ \pr\rond\partial\alpha=\partial(\pr\rond\alpha)=p_2^*(\pr\rond\alpha)
-m^*(\pr\rond\alpha)+p_1^*(\pr\rond\alpha) .\]
From $\pr\rond\Ad_g=\pr$ and $R_g^*\alpha=\Ad_{g\inv}\alpha$, it follows that 
$R_g^*(\pr\rond\alpha)=\pr\rond\alpha$. This, together with the relations 
$p_2\rond R_{(g,h)}=R_h\rond p_2$ and $p_1\rond R_{(g,h)}=R_g\rond p_1$ implies that $p_2^*(\pr\rond\alpha)$ and $p_1^*(\pr\rond\alpha)$ are invariant under the $G\times G$-action $\tgm_2$. 
Given a smooth path $t\mapsto \big(x_t,y_t\big)$ in $\tgm_2$, one also gets 
\begin{align*}
& R_{(g,h)}^* m^* (\pr\rond\alpha) \Big( \tfrac{d}{dt} \big(x_t,y_t\big) \big|_0 \Big) \\
=& (\pr\rond\alpha) \Big( \tfrac{d}{dt} x_t g y_t h \big|_0 \Big) \\ 
=& (\pr\rond\alpha) \Big( \tfrac{d}{dt} x_t y_t g^{r(y_t)} h \big|_0 \Big) \\ 
=& (\pr\rond\alpha) \Big( \tfrac{d}{dt} x_t y_t g^{r(y_0)} h \big|_0 \Big) + (\pr\rond\alpha) \Big( \tfrac{d}{dt} x_0 y_0 g^{r(y_t)} h \big|_0 \Big) 
\end{align*}
While the first term of the R.H.S. is equal to $m^*(\pr\rond\alpha) \Big( \tfrac{d}{dt} \big(x_t,y_t\big) \big|_0 \Big)$ since $(\pr\rond\alpha)$ is $G$-invariant, the second term vanishes. 
Indeed, using $\alpha(\hat{\xi})=\xi$, $R_{h}^*\alpha=\Ad_{h\inv}\rond\alpha$ and $\pr\rond\Ad_g=\pr$, we obtain that 
\begin{multline*}
(\pr\rond\alpha) \Big( \tfrac{d}{dt} x_0 y_0 g^{r(y_t)} h \big|_0 \Big) 
= \pr \Big( \tfrac{d}{dt} g^{r(y_t)} (g^{r(y_0)})\inv \big|_0 \Big) \\
= \pr \Big( \tfrac{d}{dt} {\big(g^{r(y_0)}\big)}^{r(y_0)\inv r(y_t)} (g^{r(y_0)})\inv \big|_0 \Big)
\end{multline*}
and the claim follows from~\eqref{iets}. 
Hence $R_{(g,h)}^* m^* (\pr\rond\alpha)=m^* (\pr\rond\alpha)$. 
Therefore, $\pr\rond\partial\alpha$ is $(G\times G)$-invariant.

One also has 
\begin{align*}
 &\pr\rond\partial\alpha\Big(\tfrac{d}{dt} (xe^{t\xi},ye^{t\eta})\big|_0\Big) \\
=& \pr\Big( \alpha\big(\tfrac{d}{dt} ye^{t\eta}\big|_0\big) 
- \alpha\big(\tfrac{d}{dt} x e^{t\xi} y e^{t\eta} \big|_0\big) 
+ \alpha\big(\tfrac{d}{dt} xe^{t\xi}\big|_0\big) \Big) \\ 
=& \pr\Big( \eta 
- \alpha\big(\tfrac{d}{dt} xy e^{t\Ad_{r(y)}\inv\xi} e^{t\eta} \big|_0\big) 
+ \xi \Big) \\ 
=& \pr(\eta-\Ad_{r(y)}\inv\xi-\eta+\xi) \\ 
=& \pr(\xi) - \pr(\Ad_{r(y)}\inv\xi) \\ 
=& 0 
.\end{align*}
Hence the 1-form $\pr\rond\partial\alpha\in\OO^1\big(\tgm_2,\mfzg\big)$ is basic 
with respect to the principal $(G\times G)$-bundle $\tgm_2\to\gm_2$.

\alinea{\ref{NewOrleans}} Clearly, one has $i_{\hat{\xi}}(\alpha_1-\alpha_2)=0$ and $R_g^*\big(\pr\rond(\alpha_1-\alpha_2)\big)=\pr\rond(\alpha_1-\alpha_2)$. 
Thus $\pr\rond(\alpha_1-\alpha_2)=\phi^*A$, where $A\in\OO^1\big(\gm;\mfzg\big)$. 
It follows that 
\begin{multline*} 
\phi^*(\Omega_{\alpha_1}-\Omega_{\alpha_2}) 
= \pr\big(d(\alpha_1-\alpha_2)+\partial(\alpha_1-\alpha_2)\big) 
= d(\phi^* A)+\partial(\phi^* A) \\ = \phi^*(dA+\partial A) 
\end{multline*}
and $\Omega_{\alpha_1}-\Omega_{\alpha_2}=dA+\partial A\in B^3\big(\gm_\bullet;\mfzg\big)$.
\end{proof}

\begin{remark}\label{R:DDdual}
The Dixmier--Douady class $\DDcent{\alpha}$ of a central $G$-extension identifies 
with a linear map $\mfzg^*\to H^3(\boldGamma)$ by composition with the canonical biduality 
homomorphism $\mfzg \to {\mfzg^{*}}^*$.
\end{remark}

\begin{remark}\label{rem:DDS1} 
When the group $G$ is abelian, then $\tgm/Z(G)=\gm$, $\tgm'=\tgm$ and the projection map 
$\mfg\xto{\pr}\mfzg=\mfg$ is the identity. In particular, when $G=S^1$, the Dixmier--Douady class given 
by Proposition~\ref{P:DDcent} coincides with the Dixmier--Douady class defined in~\cite{BehrendXu}.
\end{remark}

\subsection{Main theorem}

Let $\tgm\to \gm \toto M$ be a central $G$-extension of Lie groupoids.

According to Proposition~\ref{P:red}, we obtain a  universal characteristic map $\Duniv{\Phi}:H^3(\Crossed{Z(G)}{1})\to H^3(\boldGamma)$. According to~\cite{GinotXu}, $H^3(\Crossed{Z(G)}{1})$ is isomorphic to $\mfzg^*$ if $G$ is compact. Thus we obtain a map $\Duniv{\Phi}:\mfzg^*\to H^3(\boldGamma)$ which, by duality, defines the {\em universal characteristic class} $\Duniv{\Phi}\in H^3(\boldGamma)\otimes \mfzg$.  

Our main theorem is 
\begin{theorem}\label{T:CC=DD}
Let $G$ be a compact connected Lie group.
For any central $G$-extension of Lie groupoids 
$\tgm\to \gm \toto M$, the universal characteristic class coincides with the Dixmier--Douady class.
\end{theorem}

\begin{remark}
The theorem above may be considered as a higher analogue of Chern--Weil theory 
where the characteristic classes of a principle bundle
can be expressed by geometric data such as connection and
curvature. In particular when $G$ is $S^1$, Theorem~\ref{T:CC=DD}
is a higher analogue of the following well known fact: 
the Chern class of an $S^1$-bundle can be computed from its curvature. 
The latter played an important role in geometric quantization of symplectic manifolds. 
We refer the interested reader to~\cite{WX91} for prequantization of symplectic groupoids, 
which may be interpreted as a construction of a central $S^1$-extension whose Dixmier--Douady class is the 
prescribed symplectic form. 
\end{remark}

As a corollary we obtain
\begin{corollary}
Let $G$ be a compact connected Lie group. The Dixmier--Douady class of any central 
$G$-extension of Lie groupoids is an integral class. 
\end{corollary}

\begin{proof}
The singular cohomology with integer coefficients $H\upcom({\boldGamma};\ZZ) $ of  a 2-groupoid $\boldGamma$ 
is defined as for  de Rham cohomology, substituting the singular cochain complex to the de Rham forms in the constructions 
of Section~\ref{S:cohomology}, see Remark~\ref{rm:coefficient}. Therefore, given any principal $[Z(G)\to 1]$-bundle $\gpbd$ 
over $\boldGamma$, we can construct an integer valued universal characteristic homomorphism $\Duniv{\gpbd}: H\upcom(\Crossed{Z(G)}{1},\mathbb{Z}) \to H\upcom({\boldGamma},\ZZ)$ as in Section~\ref{S:characteristicmap}, 
see Remark~\ref{rm:CCcoefficient}. According to the computations in~\cite{GinotXu}, the image of $H^3(\Crossed{Z(G)}{1},\mathbb{Z})$ under the canonical morphism $H^3(\Crossed{Z(G)}{1},\mathbb{Z})\to H^3(\Crossed{Z(G)}{1})$ is the lattice in $\mfzg^*$ generated by the fundamental classes of each circle component of $Z(G)\cong S^1\times\cdots\times S^1$. 
Now the result follows from Theorem~\ref{T:CC=DD}. 
\end{proof}

\subsection{The case of central $S^1$-extensions}

In this section, we establish Theorem~\ref{T:CC=DD} in the case $G=S^1$. 

Assume $\tgm\xto{\phi}\Gamma\toto M$ is a central $S^1$-extension. 
We consider the following four 2-groupoids: 
\begin{align*}
& \boldA :\quad \Gamma\toto\Gamma\toto M 
&& \boldB :\quad \tgm\times_{\Gamma}\tgm
\toto\tgm\toto M \\ 
& \boldC :\quad Z(G)\toto\ops\toto\ops 
&& \boldD :\quad \tgm\toto\tgm\toto M
\end{align*}
The central extension $\phi$ determines
the (generalized) morphisms 
\[ \xymatrix{ \boldD \ar[r]^{\phi} & \boldA } \qquad \text{and} \qquad 
\xymatrix{ \boldA & \boldB \ar[l]_{\phi}^{\sim} \ar[r]^f & \boldC } \] 
At the nerve level, we get 
\[ \xymatrix{ 
N_2(\boldD) \ar@<6pt>[d]^{\dD_2}\ar@<0pt>[d]|{\dD_1}\ar@<-6pt>[d]_{\dD_0} \ar[r]^{\phi} & 
N_2(\boldA) \ar@<6pt>[d]^{\dA_2}\ar@<0pt>[d]|{\dA_1}\ar@<-6pt>[d]_{\dA_0} & 
N_2(\boldB) \ar@<6pt>[d]^{\dB_2}\ar@<0pt>[d]|{\dB_1}\ar@<-6pt>[d]_{\dB_0} 
\ar[l]_{\phi} \ar[r]^{f} &  
N_2(\boldC) \ar@<6pt>[d]^{\dC_2}\ar@<0pt>[d]|{\dC_1}\ar@<-6pt>[d]_{\dC_0} \\ 
\tgm \ar[r]^{\phi} & 
\Gamma & 
\tgm \ar[l]_{\phi} \ar[r]^{f} &  
\ops 
} \]
where, according to Remark~\ref{paramN2}, 
\begin{gather*} 
N_2(\boldA)=\Gamma\times_{t,\Gamma_0,s}\Gamma, \qquad  
N_2(\boldC)=Z(G), \qquad 
N_2(\boldD)=\tgm\times_{t,\Gamma_0,s}\tgm, \\ 
N_2(\boldB)=\genrel{(a,b,c)\in\tgm^3}{\phi(a)
=\phi(b)\text{ and }s(a)=s(b)=t(c)} \\  
\phi:N_2(\boldB)\to N_2(\boldA):(a,b,c)\mapsto \big(\phi(a),\phi(c)\big) \\ 
f:N_2(\boldB)\to N_2(\boldC):(a,b,c)\mapsto ab\inv 
\end{gather*}
and the face maps are given by 
\begin{align*}
\dA_0(a,c)&=a & \dA_1(a,c)&=ac & \dA_2(a,c)&=c \\ 
\dB_0(a,b,c)&=a & \dB_1(a,b,c)&=bc & \dB_2(a,b,c)&=c \\ 
\dD_0(a,c)&=a & \dD_1(a,c)&=ac & \dB_2(a,c)&=c 
\end{align*} 

We will need one more map: 
\[ \lion:N_2(\boldB)\to N_2(\boldD):(a,b,c)\mapsto(a,c) \]

\begin{lemma} One has 
\begin{gather} \dB_0=\dD_0\rond\lion, \qquad \dB_2=\dD_2\rond\lion, \label{un} \\ 
\text{and} \qquad \dD_1\rond\lion(a,b,c)=f(a,b,c)\cdot\dB_1(a,b,c), 
\qquad \forall (a,b,c)\in N_2(\boldB). \label{deux} 
\end{gather}
\end{lemma}

\begin{lemma}\label{subtractpartials}
For any pseudo-connection $\theta\in\OO(\tgm)$ on the central $S^1$-extension 
$\tgm\xto{\phi}\Gamma\toto\Gamma_0$, one has 
\[ \partial^{\boldB}\theta+f^*(dt)=\lion^*(\partial^{\boldD}\theta) \] 
Here $dt$ denotes the Maurer-Cartan (or angular) form on $S^1$.
\end{lemma}

\begin{proof}
Since $\theta\big(\frac{d}{dt}\gammatilde\cdot e^{it}\big|_0\big)=1$ and 
$\theta$ is $S^1$-invariant, it follows from~\eqref{deux} that 
\begin{equation}\label{trois} (\dD_1\rond\lion)^*\theta=(\dB_1)^*\theta+f^*dt \end{equation}
Therefore, 
\begin{align*} 
\partialB\theta-\lion^*(\partialD\theta) 
=& (\dB_1)^*\theta-\lion^*(\dD_1)^*\theta && \text{by~\eqref{un},} \\ 
=& -f^*dt && \text{by~\eqref{trois}.} 
\qedhere \end{align*}
\end{proof}

According to Proposition~\ref{P:DDcent}, the connection $\theta$ induces a cocycle $\Omega_\theta \in \ZDR{3}(\boldA)$.
\begin{theorem}\label{DDareequal}
$\phi^*[\Omega_{\theta}]=f^*[dt]$ in $H^3(\boldB)$
\end{theorem}

\begin{proof}
By construction, the cocycle $\Omega_{\theta}$ is the sum $\Omega_{\theta}=\eta+\omega$, where $\phi^*(\eta)=\partialD\theta$ and $\phi^*(\omega)=\ddR\theta$.
\begin{align*}
\partialB\theta+d\theta 
=& \lion^*(\partialD\theta)-f^*(dt)+d\theta 
&& \text{by Lemma~\ref{subtractpartials}} \\ 
=& \lion^*(\phi^*\eta)-f^*(dt)+\phi^*\omega 
&& \\ 
=& \phi^*(\eta+\omega)-f^*(dt) && 
\qedhere \end{align*} 
\end{proof}

Now, Theorem~\ref{T:CC=DD} in the case $G=S^1$ follows from Theorem~\ref{DDareequal} since $\Duniv{\phi}=\big(\phi^*\big)^{-1}(f^*(dt))$ and $\DDcent{\phi}=[\Omega_{\theta}]$ (by Proposition~\ref{P:DDcent}). 

\subsection{Proof of Theorem~\ref{T:CC=DD}}
By \cite{NADG} (see also Section~\ref{S:DDcent}), the central $G$-extension of Lie groupoids $\tgm\to \gm \toto M$ induces a central $Z(G)$-extension $\tgm'\to \gm\toto M$, where $\tgm'=q^{-1}(\sigma(\gm))$. We recover $\tgm$ from $\tgm'$ by the formula $\tgm\cong\frac{\tgm'\times G}{Z(G)}$, where the  $Z(G)$-action on $\tgm'\times G$ is given by $(x,g)\cdot z = (x\cdot z^{-1}, z\cdot g)$ for $x\in \tgm'$, $g\in G$ and $z\in Z(G)$. The natural inclusion $\tau: \tgm'\to \tgm$ coincides with the map $x\mapsto [x,1_G] \in\frac{\tgm'\times G}{Z(G)}$. Since $G$ is compact, $\mfzg$ is reductive and we have the Lie algebra morphism $\pr: \mfg\to \mfzg$.
\begin{lemma}\label{L:pralpha}
Let $\alpha \in \Omega^1(\tgm,\mfg)$ be a connection 1-form on the right principal $G$-bundle $\tgm \to \gm$. Then $\alpha':=\pr(\tau^*(\alpha))\in \Omega^1(\tgm',\mfzg)$ is a connection 1-form for the right principal $Z(G)$-bundle $\tgm' \to \gm$.
\end{lemma}

\begin{proof}
Since the inclusion $\tau:\tgm'\into\tgm$ is $Z(G)$-equivariant, we have 
\[ \alpha'(\hat{\eta}_x)=\pr\rond\alpha\rond\tau_*(\hat{\eta}_x)=
\pr\rond\alpha\rond(\hat{\eta}_{\tau(x)})=\pr(\eta)=\eta \] 
for all $\eta\in\mfzg$ and $x\in\tgm'$.
Similarly, for any $h\in Z(G)$, we have 
\begin{equation*}
R_h^*(\alpha')=R_h^* (\pr(\tau^*(\alpha)))
=\pr \tau^* R_h^*\alpha=
\pr \tau^* Ad_{h\inv}\alpha=\alpha'.
\qedhere \end{equation*}
\end{proof}

Since $\tgm' \to \gm \toto M$ is a $Z(G)$-central extension,
 by Lemma~\ref{L:pralpha} and  Proposition~\ref{P:DDcent}, we have the Dixmier--Douady class $\DDcent{\alpha'}\in H^3(\gm)\otimes \mfzg$.
\begin{proposition}\label{P:DDZGtoG} 
We have $\DDcent{\alpha'}=\DDcent{\alpha}$.
\end{proposition}

\begin{proof} By Lemma~\ref{L:pralpha} and Proposition~\ref{P:DDcent}.(b), we can use the 1-form $\alpha'$ to calculate the Dixmier--Douady class of $\tgm' \to \gm \toto M$. By construction we have a commutative diagram of groupoid morphisms
\[ \xymatrix{\tgm' \ar[rd]^{\phi'} \ar[d]_{\tau} & \\
\tgm \ar[r]^{\phi} & \gm  .} \] 

 According to Proposition~\ref{P:DDcent}.(a), the Dixmier--Douady class $\DDcent{\alpha}$ is the cohomology class of the cocycle $ [\Omega_{\alpha}]$  defined by the identity
\begin{equation}\label{eq:DDalpha}
\pr \big(d\alpha + \partial \alpha\big) =\phi^*(\Omega_\alpha).
\end{equation}
Applying $\tau^*$ to Equation~\eqref{eq:DDalpha}, we get
\begin{align*}
\tau^*\pr \big(d\alpha + \partial \alpha\big) &= \tau^*\phi^*(\Omega_\alpha) \\
\pr \big(d +\partial\big) \tau^* \alpha &= {\phi'}^*(\Omega_\alpha) \\
\big(d +\partial\big)\pr \tau^* \alpha &= {\phi'}^*(\Omega_\alpha).
\end{align*}
Therefore, by Proposition~\ref{P:DDcent}, $\DDcent{\alpha'}=[\Omega_{\alpha}]=\DDcent{\alpha}$.
\end{proof}

Since $G$ is compact its center is the quotient $(Z_0(G)\times C)/N$, where $Z_0(G)$ is the connected component of
$1_G$ in $Z(G)$ and $C$, $N$ are finite. We fix an isomorphism of Lie groups $Z_0(G)\cong S^1\times \cdots \times S^1$
(with $n$-factors). We thus obtain isomorphisms 
$\mfzg\cong \RR e_1\oplus \cdots \oplus \RR e_n$ and $H^3(\Crossed{Z(G)}{1})\cong \RR dt_1 \oplus \cdots \oplus \RR dt_n$. 
Let $\pr_i: \mfzg \to \RR e_i$ ($i=1\dots n$) be the natural projection.

\begin{lemma}\label{L:priCC} We have 
$\pr_i\big(\Duniv{\Phi} \big)=\Duniv{\Phi}(dt_i)$ in $H^3(\gm)$.
\end{lemma}

\begin{proof} Let $(\xi_1,\dots, \xi_n)$ be the dual basis of $(e_1,\dots,e_n)$ in $\mfzg^*$.
According to~\cite{GinotXu}, the generator $dt_i$ is  the left invariant vector field 
$\xi_i^L\in \Omega^1(Z(G)) \subset \Omega^3(\Crossed{Z(G)}{1})$ associated to $\xi_i$.  The Lemma follows.
\end{proof}

\begin{proposition}\label{P:CCtoDDZG} 
We have $\Duniv{\Phi}=\DDcent{\alpha'}$
\end{proposition}

\begin{proof}
By linearity and Lemma~\ref{L:pralpha}, it is sufficient to prove that for all $i=1\dots n$, one has
\begin{equation}\label{eq:pri}
\pr_i\big(\DDcent{\theta} \big)= \pr_i\big(\Duniv{\Phi} \big) = \Duniv{\Phi}(dt_i) \quad \text{ (by Lemma~\ref{L:priCC})}.
\end{equation}
The proof of Equation~\eqref{eq:pri} is similar to that of Theorem~\ref{DDareequal}.  
\end{proof}

\begin{proof}[Proof of Theorem~\ref{T:CC=DD}]
By Proposition~\ref{P:DDZGtoG} and Proposition~\ref{P:CCtoDDZG} we obtain 
\[ \Duniv{\Phi}=\DDcent{\alpha'}=\DDcent{\alpha} \]
and Theorem~\ref{T:CC=DD} follows.
\end{proof}

\begin{example}
Let $G$ be a simple compact Lie group and let $LG$ denote its loop group. 
Suppose that the Lie algebra $\mathfrak{g}$ of $G$ is endowed with $m$ invariant 
non-degenerate bilinear symmetric forms ${\langle - , - \rangle}_i$ ($i=1,\dots,m$) 
and assume that the Lie algebra 2-cocycle $\beta\in\Lambda^2(L\mathfrak{g}^*)\otimes\mathbb{R}^m$ defined by 
\[ \beta(X,Y)=\Big(\frac{1}{2\pi}\int_0^{2\pi}{\langle X(s),Y'(s)\rangle}_1\,ds,\cdots,
\frac{1}{2\pi}\int_0^{2\pi}{\langle X(s), Y'(s)\rangle}_m \,ds \Big) \] 
is integral on every factor (i.e.\ the associated closed 2-form is). 
It thus gives rise to a central extension 
 \[ T^m \to \widetilde{LG} \to LG \] of the loop group by a torus $T^m$ of dimension $m$.
 The extension being central, the adjoint action of $\widetilde{LG}$ on its Lie algebra 
 $\widetilde{L\mathfrak{g}}=L\mathfrak{g}\oplus \mathbb{R}^m$ descends to an action 
 on $L\mathfrak{g}$ (which is compatible with the adjoint action of $LG$ on $L\mathfrak{g}$). 
 Hence we have a central $T^m$-extension
 \begin{equation}\label{eq:exTmGerbe} \big(\widetilde{LG} \times L\mathfrak{g} \toto  L\mathfrak{g} \big) 
\longrightarrow \big(LG  \times L\mathfrak{g} 
\toto  L\mathfrak{g}\big)
 \end{equation}of the associated transformation Lie groupoids.
 Since $\big(LG  \times L\mathfrak{g} 
\toto  L\mathfrak{g}\big)$ is Morita equivalent\footnote{The equivalence is induced by the projection $(\ev_0,\Hol)$, 
where $\ev_0$ is the evaluation map which takes a loop in $G$ to its value at $0$ and 
$\Hol: L\mathfrak{g}\to G$ is the holonomy.} to the transformation groupoid $G\times G \toto G$
(where $G$ acts on itself by conjugation),
the extension~\eqref{eq:exTmGerbe} defines a central $T^m$-gerbe on $G\times G \toto G$. 

We compute the universal characteristic class of this central gerbe using the Dixmier--Douady class as follows. Note that the coadjoint action 
of $\widetilde{LG}$ on $\widetilde{L\mathfrak{g}}^*$ also restricts to $L\mathfrak{g}^*$ (identified with the affine hyperplane 
$\{(x,1,\dots,1); x\in L\mathfrak{g}\}$).
Let $\alpha=\alpha_1\oplus \cdots \oplus \alpha_m \in \Omega^1(\widetilde{LG}\times L\mathfrak{g})\otimes \mathbb{R}^m$ 
be the direct sum of the  dual (using on each coordinate of 
$\mathbb{R}^m$ the identification of $L\mathfrak{g}$ with $L\mathfrak{g}^*$ given by the form $\langle \, , \rangle_i$ )
of the restriction of the Liouville $1$-form on $\widetilde{LG}\times \widetilde{L\mathfrak{g}}^*$ to 
$\Omega^1(\widetilde{LG}\times L\mathfrak{g}^*)$. 
By Proposition~\ref{P:DDcent}  and Lemma~\ref{L:priCC}, we obtain that the Dixmier--Douady class of the gerbe  is a sum 
$\DDcent{\alpha}= [\Omega_{\alpha_1}]\oplus \cdots \oplus [\Omega_{\alpha_m}]$. By computation 
(see~\cite{BehrendXuZhang}*{Proposition~3.2}), we see
that each class $[\Omega_{\alpha_i}]$ is equal in $H^3(G\times G\toto G)$ to the class of $[b_i+\lambda_i]$ where 
$\lambda_i$ is the bi-invariant $3$-form coresponding to $\frac{1}{12} {\langle -, [-,-]\rangle}_i\in \Lambda^3 \mathfrak{g}^*$ 
and 
\begin{multline*} 
b_i= -\frac{1}{2}[{\langle Ad_{x} g^*(\theta_{MC}), g^*(\theta_{MC})\rangle}_i 
+ {\langle g^*(\theta_{MC}, x^*(\theta_{MC}+\theta_{MC}^r) \rangle}_i] \\ 
\in \Omega^2(G\times G)\subset \Omega^3(G\times G\toto G) 
,\end{multline*}
where $\theta_{MC}$ and $\theta_{MC}^r$ are respectively the left and right Maurer-Cartan forms 
and $(g,x)$ denotes the coordinates on $G\times G$.

Since $G$ is a simple compact group, each ${\langle -, - \rangle}_i$ is an integer multiple 
\[ \langle -, -\rangle_i = a_i \langle -, -\rangle_{\bas} \] of the basic form of $G$. 
Hence $[b_i+\lambda_i]=a_i[b_{\bas}+\lambda_{\bas}]$ (where $b_{\bas},\lambda_{\bas}$ are defined as 
above). Further, $[b_{\bas}+\lambda_{\bas}]$ is precisely an integral generator of $H^3(G\times G\toto G)$.

It follows that the characteristic class of the gerbe is 
\[ \DDcent{\alpha} = a_1\oplus \cdots \oplus a_m \in \mathbb{R}^m .\]
\end{example}

\section*{Acknowledgments} The authors are grateful to Andr\'e Haefliger, Jim Stasheff, Urs Schreiber, Danny Stevenson, Ping Xu and the referee for many useful discussions and suggestions and are beholden to Damien Broka for contributing improvements to Section~\ref{S:gpbdtoGextension}.

\end{document}